\documentclass[11pt,a4paper]{article}
\usepackage[numbers,sort&compress]{natbib}
\usepackage[colorlinks,
            linkcolor=blue,
            anchorcolor=green,
            citecolor=magenta
            ]{hyperref}
\usepackage{mathrsfs,amssymb,amsfonts}
\usepackage{graphicx}
\usepackage{amsmath,amssymb,latexsym,color}
\usepackage[mathscr]{eucal}
\usepackage[numbers]{natbib}
\usepackage{bm}
\usepackage[top=1in, bottom=1in, left=1.25in, right=1.25in]{geometry}
\makeatletter

\newcommand{\Rmnum}[1]{\expandafter\@slowromancap\romannumeral #1@}
\makeatother
\newtheorem{thm}{Theorem}[section]

\newtheorem{lem}[thm]{Lemma}

\newtheorem{example}{Example}
\newtheorem{remark}{Remark}[section]

\newcommand{\qed}{\hfill\Box\medskip}
\usepackage{CJK}
\usepackage{geometry}
\begin{document}
\begin{CJK*}{GBK}{song}
 \setlength{\baselineskip}{11.5pt}
\renewcommand{\abovewithdelims}[2]{
\genfrac{[}{]}{0pt}{}{#1}{#2}}
%%%%%%%%%%%%%%%%%%%%%%%%%%%%%%%%%%%%%%%%%%%%%%%%%%%%%%%%%%%%%%%%%%%%%%%%%%%%%%%%%%%%%%%%
%%%%%%%%%%%%%%%%%%%%%%%%%%%%%%%%%%%%%%%%%%%%%%%%%%%%%%%%%%%%%%%%%%%%%%%%%%%%%%%%%%%%%%%%

\title{\bf The super spanning connectivity of arrangement graphs}

\author{Pingshan Li \quad  Min Xu\footnote{Corresponding author. \newline {\em E-mail address:} xum@mail.bnu.edu.cn (M. Xu) .}\\
{\footnotesize   \em  Sch. Math. Sci. {\rm \&} Lab. Math. Com. Sys., Beijing Normal University, Beijing, 100875,  China} }
 \date{}
 \date{}
 \maketitle

\begin{abstract}
A $k$-container  $C(u, v)$ of a graph $G$ is a set of $k$ internally disjoint paths between $u$ and $v$. A $k$-container $C(u, v)$ of $G$ is
 a $k^*$-container if it is a spanning subgraph of $G$. A graph $G$ is $k^*$-connected if there exists a $k^*$-container between
 any two different vertices of G.  A $k$-regular graph $G$ is super spanning connected if $G$ is $i^*$-container for all $1\le
  i\le k$. In this paper, we prove that the arrangement graph $A_{n, k}$ is super spanning connected if $n\ge 4$ and $n-k\ge 2$.

\medskip
\noindent {\em Key words:} Hamiltonian; Hamiltonian connected; $k^*$-connected; Arrangement graphs.

\medskip
%\noindent {\em 2010 MSC:} 05C25; 05C15.
\end{abstract}

\section{Introduction}
In the field of parallel and distributed systems, interconnection networks are an important research area. Typically, the topology
of a network can be represented as a graph in which the vertices represent processors and the edges represent communication links.

For graph definitions and notations, we follow \cite{J.A.Bondy}.
A graph $G$ consists of a vertex set $V(G)$ and an edge set $E(G)$, where an edge is an unordered pair of distinct vertices of $G$.
A path  $P$ of length $k$ from $x$ to $y$ is a finite sequence of distinct vertices $\langle v_0, v_1, \cdots, v_k\rangle$, such that $x=v_0, y=v_k$, and $(v_i, v_{i+1})\in E$ for $0\le i\le k-1$. We also represent path $P$ as $\langle v_0, v_1, \cdots, v_i, Q, v_j, v_{j+1}, \cdots, v_k\rangle$, where $Q$ is the path $\langle v_i, v_{i+1}, \cdots, v_j\rangle$.  In particular, if $i=j$, we can still represent the path as $\langle v_0, v_1, \cdots, v_i, Q, v_i, v_{j+1}, \cdots, v_k\rangle$.

A spanning subgraph of $G$ is a subgraph with vertex set $V(G)$. A Hamiltonian graph is a graph with a spanning cycle.
A graph is Hamiltonian connected if there exists a spanning path joining any two different vertices.

A $k$-container $C(u, v)$ of $G$ is a set of $k$ internally disjoint  paths between $u$ and $v$, and
the connectivity of $G$, $\kappa(G)$, is the minimum size of a vertex set  $S$ such that $G-S$ is disconnected or has only one vertex.
It follows from Menger's Theorem \cite{K.Menger} that there is a $k$-container between any two distinct vertices of $G$ if $G$ is $k$-connected.

A $k$-container $C(u, v)$ of $G$ is a $k^*$-container  if it is a spanning subgraph of $G$. A graph is $k^*$-connected if there exists a $k^*$-container between any two distinct vertices.
By this definition, the concept of Hamiltonian connected is the same as $1^*$-connected and the concept of Hamiltonian is the same as $2^*$-connected. Thus, the concept of $k^*$-connected is a hybrid concept of connectivity and Hamiltonicity.
The study of $k^*$-connected graphs is motivated by the globally $3^*$-connected graphs proposed by M. Albert et al. \cite{M.Albert}.

The star graph ($S_n$ for short), which was proposed by Akers et al.\cite{S.B.Akers1}, is a well known interconnection
network. The arrangement graph\cite{Day and Tripathi}, denoted by $A_{n,k}$ , refers to a generalized version of $S_n$.
 Further, $A_{n, n-1}$ is isomorphic to the $n$-dimensional star graph $S_n$\cite{S.B.Akers} and $A_{n, 1}$ is isomorphic to the complete graph $K_n$. The arrangement graph preserves many attractive properties of $S_n$ such as the hierarchical structure, vertex and edge symmetry, simple and optimal routing, and many fault tolerance properties \cite{Day and Tripathi}. Some basic properties of $A_{n, k}$ such as average distance\cite{on}, Hamiltonicity\cite{Fault}, and embedding\cite{Xumin,KDay} have recently been computed or derived.

A graph $G$ is super spanning connected if it is $k^*$-connected for all $1\le k\le \kappa(G)$.
There are many desirable results about super spanning connected of some interconnection networks such as recursive circulant graphs\cite{C.H.Tsai},
 pancake graphs\cite{C.K.Lin},  hypercube-like network\cite{Cheng-Kuan Lina},  $(n, k)$-star graphs\cite{H.C.Hsu2}, $k$-ary $n$-cubes\cite{Y.K. Shih} and multi-dimensional tori\cite{Lijing}. 
 Since  $A_{n, n-1}\cong S_n$ is a bipartite graph with the same number of vertices in each partite set, there is no Hamiltonian path joining any two different vertices in the same part. Hence, $A_{n, n-1}$ is not super spanning connected if $n>3$.
 Therefore, we consider an arrangement graph with $n-k\ge 2$.
In this paper, we aim to prove that arrangement graphs are super spanning connected if $n\ge 4$ and $n-k\ge2$.

The rest of this paper is organized as follows. In Section 2, we introduce   arrangement graphs and discuss some of their properties. In Section 3, we prove that arrangement graphs  are super spanning connected for $n\ge 4$ and $n-k\ge2$.

\section{Arrangement graphs}

\begin{figure}[!htbp]
 \centering
 \includegraphics[width=0.35\textwidth]{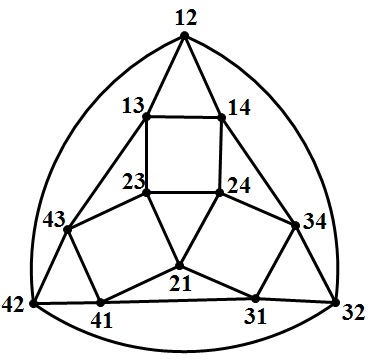}
 \caption{The arrangement graph $A_{4, 2}$}\label{a42}
 \end{figure}

Throughout this paper, we assume that $n$ and $k$ are positive integers with $n>k$. We use $\langle n\rangle$ to denote the set $\{1, 2, \cdots, n\}$. The arrangement graph $A_{n, k}$
 is a graph that has the vertex set $V(A_{n, k})=\{u=u_1u_2\cdots u_k\mid u_i\in \langle n\rangle, u_i\not=u_j$ if $i\not=j\}$ and the edge set $E(A_{n, k})=\{(p, q)\mid p, q\in V(A_{n, k})$ and $p, q$ differ in exactly one position$\}$. From the definition , we know that $A_{n, k}$ is
a regular graph of degree $k(n-k)$ with $\frac{n!}{(n-k)!}$ vertices.
 Figure \ref{a42} illustrates the arrangement graph $A_{4, 2}$.

Let $u=u_1u_2\cdots u_k\in V(A_{n, k})$. We denote $(u)_i=u_i$ as the ith coordinate of $u$
 for $1\le i\le k$. Let $v$ be a neighbor of $u$, we denote $v$ as $u^{s(u_i, x)}$ if $v=u_1u_2\cdots u_{i-1}xu_{i+1}\cdots u_k$ for $x\in \langle n\rangle\setminus\{(u)_i: i=1, 2, \cdots, k\}$.

 For $i, j\in\langle n\rangle, l\in\langle k\rangle$ and $i\not=j$, suppose that $A^{(l, i)}_{n, k}$ denotes  the subgraph of $A_{n, k}$ that is induced by $V(A^{(l, i)}_{n, k})=\{p\mid p=p_1p_2\cdots p_k$ and $p_l=i\}$. Obviously,
 $\{V(A^{(l, i)}_{n, k})\mid 1\le i\le n\}$ forms a partition of $V(A_{n, k})$ and each $A^{(l, i)}_{n, k}$ is isomorphic to $A_{n-1, k-1}$. As a result, $A_{n, k}$ can be recursively constructed from $n$ copies of $A_{n-1, k-1}$.
 We use $E^{l=i, j}$ to denote the set of edges between
 $A^{(l, i)}_{n, k}$ and $A^{(l, j)}_{n, k}$; accordingly, $E^{l=i, j}=\frac{(n-2)!}{(n-k-1)!}$. We also use $A^{(l, I)}_{n, k}$ to denote the subgraph of $A_{n, k}$ that is induced by $\cup_{i\in I}V(A^{(l, i)}_{n, k})$.

Following are some known properties about arrangement graphs.

\begin{lem}{\rm (\cite{H.C.Hsu})}\label{L1}
  The arrangement graph $A_{n,k}$ is $(k(n-k)-2)$-fault-tolerant Hamiltonian, and $(k(n-k)-3)$-fault-tolerant Hamiltonian-connected for $n-k \ge 2$.
\end{lem}

\begin{lem}{\rm (\cite{Lo and Chen})}\label{LL1}
 The arrangement graph $A_{n, k}$ is $(k(n-k)-2)$-edge-fault-tolerant Hamiltonian connected if not all faulty edges are adjacent to the same vertex.
\end{lem}

In the following, we discuss some properties that will be used in the proof of the main results.
\begin{lem}{\rm }\label{L5}
$A_{4, 2}$ is super spanning connected.
\end{lem}

\noindent{\bf Proof: } By Lemma \ref{L1}, $A_{4,2}$ is $1^*$-connected and $2^*$-connected.
We need to construct a $3^*$-container and a $4^*$-container joining any two different vertices $u$ and $v$ of $A_{4, 2}$.
Since $A_{4, 2}$ is vertex and edge transitive, without loss of generality,
we can assume that $u=12$ and $v=13, 34, 23, 21$. We list such $3^*$-containers  as follows:
\begin{center}
 \begin{tabular}{|l|}
  \hline
  $\langle 12, 13\rangle ~~~~~~~~~~~~~~$
$ \langle 12, 14, 13\rangle~~~~~~~~~~~$
 $\langle 12, 42, 43, 41, 21, 31, 32, 34, 24, 23, 13\rangle $\\
 \hline
 $\langle 12, 14, 34\rangle~~~~~~~~~~$
 $\langle 12, 42, 32, 34\rangle~~~~~~$
  $\langle 12, 13, 43, 23, 24, 21, 41, 31, 34\rangle$ \\
 \hline
  $\langle 12, 13, 43, 23  \rangle~~~~~~$
  $\langle 12, 14, 34, 24 , 23 \rangle~~$
 $\langle 12, 42, 32, 31, 41 , 21, 23 \rangle$\\
  \hline
 $\langle 12, 13, 43, 23 , 21 \rangle~~$
  $\langle 12, 14, 34, 24 , 21 \rangle~~$
 $\langle 12, 42, 32, 31, 41 , 21 \rangle$\\
  \hline

\end{tabular}
\end{center}
and such  $4^*$-containers as follows:
 %$~~~~~~~~~~~~~~~~~~~~~~~~~~$
\begin{center}
\begin{tabular}{|ll|}
  \hline
  $\langle 12, 13\rangle$&
  $\langle 12, 14 ,13\rangle$\\
$\langle 12, 42, 43 ,13\rangle$ &
 $\langle 12, 32, 34, 31, 41, 21, 24, 23, 13\rangle$ \\
 \hline
  $\langle 12, 13, 43, 23, 24, 34\rangle$&
    $\langle 12, 14, 34\rangle$\\
     $\langle 12, 32, 34\rangle$&
      $\langle 12, 42, 41, 21$, $31, 34\rangle$ \\
      \hline
       $\langle 12, 13,  23  \rangle$&
   $\langle 12, 14,  24 , 23 \rangle$\\
    $\langle 12, 32, 34 , 31, 21, 23 \rangle$&
     $\langle 12, 42, 41 , 43, 23 \rangle$ \\
     \hline
  $\langle 12, 13, 43, 23 , 21 \rangle$&
   $\langle 12, 14, 34, 24 , 21 \rangle$\\
    $\langle 12, 32, 31 , 21 \rangle$&
     $\langle 12, 42, 41 , 21 \rangle$ \\
  \hline
\end{tabular}
\end{center}

 $\qed$

\begin{lem}{\rm }\label{L2}
  Suppose that $n\ge 5, k\ge 2 $ and $n-k\ge 2$. Let $I=\{i_1, i_2, \cdots, i_m\}\subseteq \langle n\rangle$,
 then $A^{(i, I)}_{n, k}$ is Hamiltonian connected where $1\le i\le k$.
\end{lem}

\noindent{\bf Proof: } Let $u$ and $v$ be any  two distinct vertices of $A^{(i, I)}_{n, k}$, we will prove that there exists a Hamiltonian path $P$ of $A^{(i, I)}_{n, k}$ joining $u$ and $v$. If $|I|=1$, by Lemma \ref{L1}, it is true. Therefore, we assume that $|I|\ge 2$ in the next proof.

\noindent{\bf Case 1:} $(u)_i\not=(v)_i$.

\begin{figure}[!htbp]
 \centering
 \includegraphics[width=0.50\textwidth]{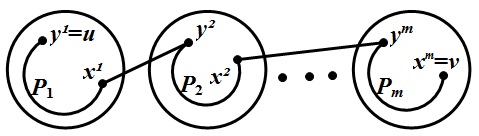}
 \caption{Illustration for case 1 of Lemma \ref{L2}}\label{L2_1}
 \end{figure}

Without loss of generality, let $(u)_i=i_1$ and $(v)_i=i_m$.
There exists at least $3$ edges between $A^{(i, i_j)}_{n, k}$ and $A^{(i, i_{j+1})}_{n, k}$ for all $1\le j\le m-1$ owning to $|E^{i=i_j, i_{j+1}}|=\frac{(n-2)!}{(n-k-1)!}\ge 3$.
An edge $(x^{i_j}, y^{i_{j+1}})\in E^{i=i_j, i_{j+1}}$ is chosen such that $(x^{i_j})_i=i_j$, $(y^{i_{j+1}})_i=i_{j+1}$ and $x^1\not=u, y^m\not=v, x^{i_j}\not=y^{i_j}$ for all $j=1, 2, \cdots, m-1$. Let $y^1=u, x^m=v$. Since $A_{n, k}$ is Hamiltonian connected , there exists a Hamiltonian
path $\langle y^{i_j}, P_j, x^{i_j}\rangle$ of $A^{(i, i_j)}_{n, k}$ joining $ y^{i_j}$ to $x^{i_j}$ for $1\le j\le m$. Hence, there exists a Hamiltonian path $P=\langle y^1, P_1, x^1, y^2, P_2, x^2, \cdots,y^m, P_m, x^m\rangle $ of $A^{(i, I)}_{n, k}$ joining $u$ to $v$. See figure \ref{L2_1} for illustration.

\noindent{\bf Case 2:} $(u)_i=(v)_i$.

\begin{figure}[!htbp]
 \centering
 \includegraphics[width=0.50\textwidth]{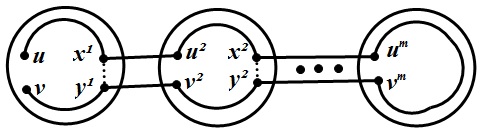}
 \caption{Illustration for case 2 of Lemma \ref{L2}}\label{L2_2}
 \end{figure}

Without loss of generality , let $(u)_i=(v)_i=i_1$.
A vertex $x^1\in V(A^{(i, i_1)}_{n, k})\setminus\{u, v\}$ is chosen such that $i_2\notin\{(x^1)_j: j=1, 2, \cdots, k\}$. Obviously,
$|\{x'\mid (x^1, x')\in E(A^{(i, i_1)}_{n, k})$ and $i_2\notin\{(x')_j: 1\le j\le k\} \}|=(n-k-1)(k-1)\ge 2 $ owing to $n\ge 5$.
Hence, there exists at least two neighbors $y^1, z^1$ of $x^1$ such that $i_2\notin\{(y^1)_j: 1\le j\le k\}\cup\{(z^1)_j:1\le j\le k\}$.
Let $e=(x^1, a)\in E(A^{(i, i_1)}_{n, k})$ and $a\notin\{y^1, z^1\}$.
By Lemma \ref{LL1}, there exists a Hamiltonian path $P_1$ of
$A^{(i, i_1)}_{n, k}-\bigcup_{x\in N_{A^{(i, i_1)}_{n, k}}(x^1)}\{(x^1, x)\}+\{(x^1, y^1),(x^1, z^1), e\}$ between  $u$ and $v$. Note that the degree of $x^1$ in $P_1$ is two,  then $\{(x^1, y^1),(x^1, z^1)\}\cap E(P_1)\ge 1$,
Without loss of generality, let $(x^1, y^1)\in E(P_1)$,
and we can represent $P_1$ as $\langle u, R_1, x^1, y^1, H_1, v\rangle$. Let $u^{j}={(x^{j-1})}^{s(i_{j-1}, i_j)}$ and $v^j={(y^{j-1})}^{s(i_{j-1}, i_j)}$ for $2\le j\le m$. Similarly, there exists a  Hamiltonian path $P_j$ of $A^{(i, i_j)}_{n, k}$ such that  $(x^j, y^j)\in E(P_j)$ and $i_{j+1}\notin \{(x^j)_l: 1\le l\le k\}\cup \{(y^j)_l: 1\le l\le k\}$ for all $2\le j\le m-1$. We can represent
$P_j$ as $\langle u^j, R_j, x^j, y^j, H_j, v\rangle$ for $1\le j\le m-1$. By Lemma \ref{L1}, there exists a Hamiltonian path $P_m$ of $A^{(i, i_m)}_{n, k}$ between $u^m$ and $v^m$. Hence, there exists a Hamiltonian path
$P=\langle u, R_1, x^1, u^2, R_2, x^2, \cdots, u^m, P_m, v^m, \cdots, y^2, H_2, v^2, y^1, H_1, v\rangle$ of $A^{(i, I)}_{n, k}$ joining $u$ to $v$. See figure \ref{L2_2} for illustration.

$\qed$

\begin{lem}{\rm }\label{L6}
For $m\in \langle n\rangle$. Suppose that  $A=\{u^1, u^2, \cdots, u^m\}, B=\{v^1, v^2, \cdots, v^m\}$, $A\cap B=\emptyset$, $A, B\subseteq V(A_{n, k})$. If there exists a number $t\in \langle k\rangle$ such that  $(u^i)_t\not=(u^j)_t$, $(v^i)_t\not=(v^j)_t$ for $1\le i\not=j\le m$,
then there exists $m$ disjoint paths $H_1, H_2, \cdots, H_m$ from $A$ to $B$ such that $V(\cup_{j=1}^{m} H_j)=V(A_{n, k})$.
\end{lem}

\begin{figure}[!htbp]
 \centering
 \includegraphics[width=0.65\textwidth]{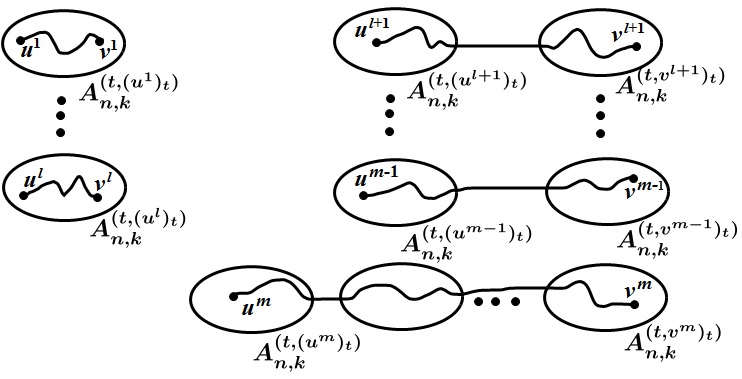}
 \caption{Illustration for  Lemma \ref{L6}}\label{L_6}
 \end{figure}

\noindent{\bf Proof: } We partite $A_{n, k}$ to $\cup_{i=1}^{n} A^{(t, i)}_{n, k}$. Suppose that $|\{(u^i)_t: 1\le i\le m\}\cap \{(v^i)_t: 1\le i\le m\}|=l$. Without loss of generality, we can assume that $(u^i)_t=(v^i)_t$ for $1\le i\le l$. By Lemma \ref{L2}, there exists a Hamiltonian path  $H_i$ of $A^{(t, (u^i)_t)}_{n, k}$ joining $u^i$ to $v^i$ for $1\le i\le l$ and a Hamiltonian path $H_j$ of $A^{(t, \{(u^j)_t, (v^j)_t\})}_{n, k}$ joining $u^j$ to $v^j$ for $l+1\le j\le m-1$.
Let $I=\langle n\rangle\setminus (\{(u^i)_t: 1\le i\le m-1\}\cup \{(v^i)_t: 1\le i\le m-1\})$. By Lemma \ref{L2}, there exists a Hamiltonian path $H_{m}$ of $A^{(i, I)}_{n, k}$ joining $u^m$ to $v^m$. Obviously, $H_1, H_2, \cdots, H_m$ form the desired paths. See figure \ref{L_6} for illustration.
 $\qed$

\section{The super spanning connectivity of arrangement graphs}

\begin{lem}{\rm }\label{TH1}
 Suppose that $n\ge 4, n-k\ge2$, then $A_{n, k}$ is $l^*$-connected for $(n-k)(k-1)+1\le l\le(n-k)k$.
\end{lem}

\noindent{\bf Proof. }We prove the Lemma  by induction.

\noindent{\bf Basis step:} Since $A_{n, 1}$ is isomorphic to the complete graph $K_n$ and by Lemma \ref{L5}, $A_{4,2}$ is super spanning connected. Thus, the result holds for $A_{n, 1}$ and $A_{4, 2}$.

\noindent{\bf Induction step:} Suppose that $A_{n-1, k-1}$ is $(n-k)(k-1)^*$-connected.

We need to find an $l^*$-container between any two different vertices $u$ and $v$ of $A_{n, k}$ with $k\ge 2$ for $(n-k)(k-1)+1\le l\le (n-k)k$.
 We use $U$ to denote the set $\{(u)_i: 1\le i\le k\}$ and $V$ to denote the set $\{(v)_i: 1\le i\le k\}$.

\noindent{\bf Case 1:} $\{i\mid (u)_i=(v)_i: 1\le i\le k\}\not=\emptyset$.

Without loss of generality, let $(u)_k=(v)_k=\alpha$.
Suppose that:
\begin{center}
$\begin{array}{rl}
  U\cap V=&\{x_1, x_2, \cdots, x_t, \alpha\},\vspace{1.0ex} \\
   U\setminus (U\cap V)=&\{u'_{t+1}, u'_{t+2}, \cdots, u'_{k-1}\}, \vspace{1.0ex} \\
  V\setminus (U\cap V)=&\{v'_{t+1}, v'_{t+2}, \cdots, v'_{k-1}\},\vspace{1.0ex} \\
    \langle n\rangle\setminus (U\cup V)=&\{w_1, w_2, \cdots, w_{n+t-2k+1}\}.
\end{array}$
\end{center}
 We partite $A_{n, k}$ to $\cup^n_{i=1}A^{(k, i)}_{n, k}$.
 By induction, there exists an $(n-k)(k-1)^*$-container $\{P_1, P_2, \cdots$, $ P_{(n-k)(k-1)}\}$ of $A^{(k, \alpha)}_{n, k}$ joining $u$ and $v$. By Lemma \ref{L2}, there exists a Hamiltonian path $R_i$ of $A^{(k, w_i)}_{n, k}$ joining $u^{s(\alpha, w_i)}$ to $v^{s(\alpha, w_i)}$ for $1\le i\le n+t-2k+1$ and a Hamiltonian path $H_j$ of $A^{(k, v'_j)}_{n, k}\cup A^{(k, u'_j)}_{n, k}$ joining $u^{s(\alpha, v'_j)}$ to $v^{s(\alpha, u'_j)}$ for $t+1\le j\le k-1$.

(a) $(n-k)(k-1)+1\le l\le (n-k)(k-1)+n+t-2k+1$. (If $n+t-2k+1=0$, then (a) does not occur. )

\begin{figure}[!htbp]
 \centering
 \includegraphics[width=0.55\textwidth]{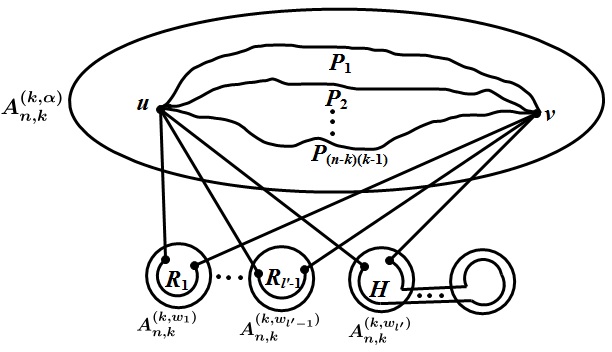}
 \caption{Illustration for case 1 (a) of Lemma \ref{TH1}}\label{TH1case1a}
 \end{figure}

Let $l=(n-k)(k-1)+l'$, then $1\le l'\le n+t-2k+1$. We set $P_{(n-k)(k-1)+i}=\langle u, u^{s(\alpha, w_i)}, R_i, v^{s(\alpha, w_i)}, v \rangle$ for $1\le i\le l'-1$. By Lemma \ref{L2}, there exists a Hamiltonian path $H$ of $A^{(k, I)}_{n, k}$ joining $u^{s(\alpha, w_{l'})}$ to $v^{s(\alpha, w_{l'})}$ where $I=\langle n\rangle-\{\alpha, w_1, w_2, \cdots, w_{l'-1}\}$. We set $P_l=\langle u, u^{s(\alpha, w_{l'})}, H$, $v^{s(\alpha, w_{l'})}, v\rangle$. Obviously, $\{P_1, P_2, \cdots, P_l\}$ forms an $l^*$-container of $A_{n, k}$ joining $u$ to $v$. See figure \ref{TH1case1a} for illustration.

(b) $(n-k)(k-1)+n+t-2k+2\le l\le (n-k)k$. (If $k-t=1$, then (b) does not occur. )

 Let $l=(n-k)(k-1)+n+t-2k+1+l'$, then $1\le l'\le k-t-1$.
 We set \vspace{1.5ex}

\begin{figure}[!htbp]
 \centering
 \includegraphics[width=0.70\textwidth]{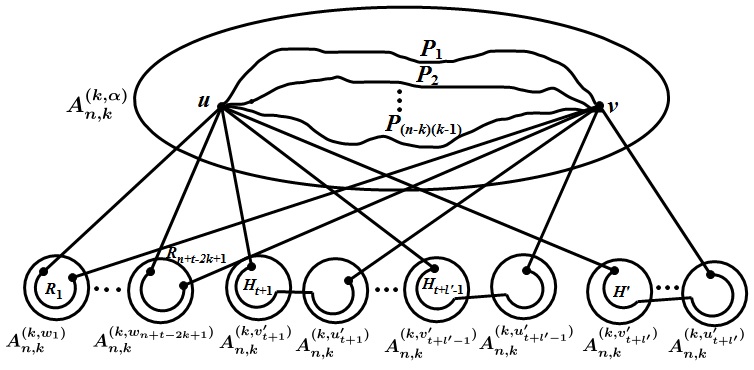}
 \caption{Illustration for case 1 (b) of Lemma \ref{TH1}}\label{TH1case1b}
 \end{figure}

$\begin{array}{rl}
  P_{(n-k)(k-1)+i}= & \langle u, u^{s(\alpha, w_i)}, R_i, v^{s(\alpha, w_i)}, v \rangle~for~ 1\le i\le n+t-2k+1 \vspace{1.0ex}, \\
  P_{(n-k)(k-1)+n-2k+1+j}= & \langle u, u^{s(\alpha, v'_j)}, H_j, v^{s(\alpha, u'_j)}, v \rangle~for~ t+1\le j\le t+l'-1.
\end{array}$\vspace{1.5ex}\\
 Let $I=\langle n\rangle\setminus\{\alpha, w_1, w_2, \cdots, w_{n+t-2k+1}, u'_{t+1}, \cdots, u'_{t+l'-1}, v'_{t+1}, \cdots, v'_{t+l'-1}\} $.
 By Lemma \ref{L2}, there exists a Hamiltonian $H'$ of $A^{(k, I)}_{n, k}$ joining $u^{s(\alpha, v'_{t+l'})}$ to
 $v^{s(\alpha, u'_{t+l'})}$. We set
 $P_l=\langle u, u^{s(\alpha, v'_{t+l'})}, H', v^{s(\alpha, u'_{t+l'})}, v \rangle$. Obviously,  $\{P_1, P_2, \cdots, P_l\}$ forms an $l^*$-container of $A_{n, k}$ joining $u$ to $v$.
 See figure \ref{TH1case1b} for illustration.

\noindent{\bf Case 2: } $\{i\mid (u)_i=(v)_i: 1\le i\le k\}=\emptyset$.

\noindent{\bf Case 2.1 : } $\{i\mid (u)_i\in V, (v)_i\in U\}\not=\emptyset$ .

Without loss of generality, we can assume that $(u)_k\in V, (v)_k\in U$. Let $(u)_k=\alpha, (v)_k=\beta$.
Suppose that:
\begin{center}
$\begin{array}{rl}
  U\cap V=&\{ x_1, \cdots, x_t, \alpha, \beta\}\vspace{1.0ex}, \\
  U\setminus (U\cap V)=&\{u'_{t+1}, u'_{t+2}, \cdots, u'_{k-2}\}\vspace{1.0ex},   \\
  V\setminus (U\cap V)=&\{v'_{t+1}, v'_{t+2}, \cdots, v'_{k-2}\}\vspace{1.0ex}, \\
  \langle n\rangle\setminus (U\cup V)=&\{w_1, w_2, \cdots, w_{n+t-2k+2}\}\vspace{1.0ex}.
\end{array}$
\end{center}
Since $n-k\ge 2$, there exists a element $\gamma\in\langle n\rangle\setminus V$.
Set $y=x_1\cdots x_tv'_{t+1}\cdots v'_{k-2}\gamma \alpha$ and $z=x_1\cdots x_tv'_{t+1}\cdots v'_{k-2}\gamma \beta$. Thus, $u\not=y, z\not=v$.
Let
$S=\langle n\rangle\setminus\{\{(y)_i: 1\le i\le k\}\cup \{\beta\}\}=\langle n\rangle\setminus\{\{(z)_i: 1\le i\le k\}\cup\{\alpha\}\}=\{s_1, s_2, \cdots, s_{n-k-1}\}$.

By induction, there exists an $(n-k)(k-1)^*$-container
$\{P_{ij}: 1\le i\le k-1, 1\le j\le n-k\}$
of $A^{(k, \alpha)}_{n, k}$ joining $u$ to $y$,  and  an $(n-k)(k-1)^*$-container $\{Q_{ij}: 1\le i\le k-1, 1\le j\le n-k\}$
of $A^{(k, \beta)}_{n, k}$ joining $z$ to $v$.
We can represent $P_{ij}$ as $\langle u, P'_{ij}, y^{ij}, y\rangle,  Q_{ij}$ as $ \langle z, z^{ij}, Q'_{ij}, v\rangle$
 for $1\le i\le k-1$ and $1\le j\le n-k$ where
 \begin{center}
$ y^{ij}=\left\{
\begin{aligned}
  & y^{s((y)_i, s_j)}: 1\le i\le k-1, 1\le j\le n-k-1, \\
 &  y^{s((y)_i, \beta)}: 1\le i\le k-1, j= n-k,
\end{aligned}
\right.$
\end{center}
\begin{center}
$ z^{ij}=\left\{
\begin{aligned}
  & z^{s((y)_i, s_j)}: 1\le i\le k-1, 1\le j\le n-k-1, \\
 &  z^{s((y)_i, \alpha)}: 1\le i\le k-1,  j= n-k.
\end{aligned}
\right.$
\end{center}

\begin{figure}[!htbp]
 \centering
  \includegraphics[width=0.88\textwidth]{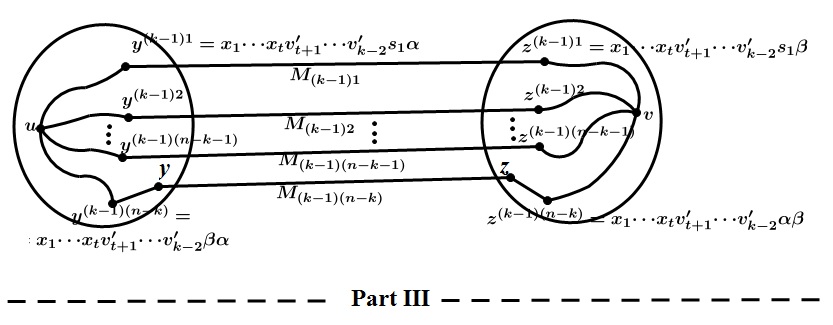}
 \caption{The $(n-k)(k-1)$  internally disjoint paths of $A_{n, k}$ of case 2.1, 2.2 of Lemma \ref{TH1}}
 \label{2_1_1}
 \end{figure}

Obviously, $(y^{ij}, z^{ij})\in E(A_{n, k})$ when  $1\le i\le k-1$ and $1\le j\le n-k-1$.
By Lemma \ref{L2}, there exists a Hamiltonian path  $R_{i}$ of $A^{(k, (y)_{i})}_{n, k}$ joining $(y^{i(n-k)})^{s(\alpha, (y)_i)}$ to $(z^{i(n-k)})^{s(\beta, (z)_i)}$ for $1\le i\le k-2$,
 As a result, there exists $(n-k)(k-1)$ internally disjoint paths $\{M_{ij}: 1\le i\le k-1, 1\le j\le n-k\}$  of $A_{n, k}$ joining $u$ to $v$ such that
\begin{center}
$V(\displaystyle\bigcup_{i=1}^{k-1}\displaystyle\bigcup_{j=1}^{n-k} M_{ij})=
V(A^{(k, \{x_1, \cdots, x_t, v'_{t+1}, \cdots, v'_{k-2}, \alpha, \beta\})}_{n, k})$
\end{center}
where \begin{center}
$ M^{ij}=\left\{
\begin{aligned}
  &\langle u, P'_{ij}, y^{ij}, z^{ij}, Q'_{ij}, v\rangle~for~1\le i\le k-1, 1\le j\le n-k-1, \\
 &\langle u, P'_{ij}, y^{ij}, (y^{ij})^{s(\alpha, (y)_i)}, R_i, (z^{ij})^{s(\beta, (z)_i)}, z^{ij}, Q'_{ij}, v\rangle~for~1\le i\le k-2, j=n-k, \\
 &\langle u, P'_{ij}, y^{ij}, y, z, z^{ij}, Q'_{ij}, v\rangle~for~i=k-1, j=n-k.
\end{aligned}
\right.$
\end{center}
 See figure \ref{2_1_1} for illustration.

\begin{figure}[!htbp]
 \centering
 \includegraphics[width=0.88\textwidth]{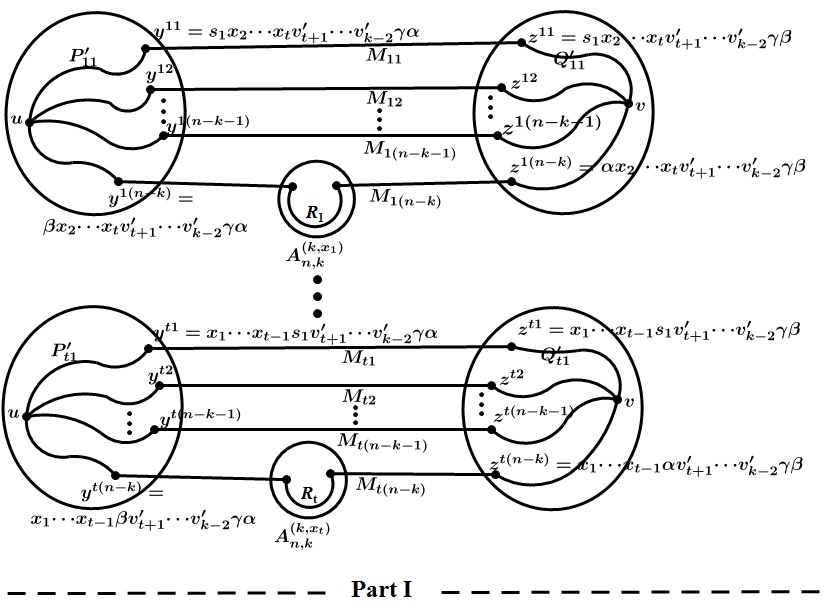}
  \includegraphics[width=0.88\textwidth]{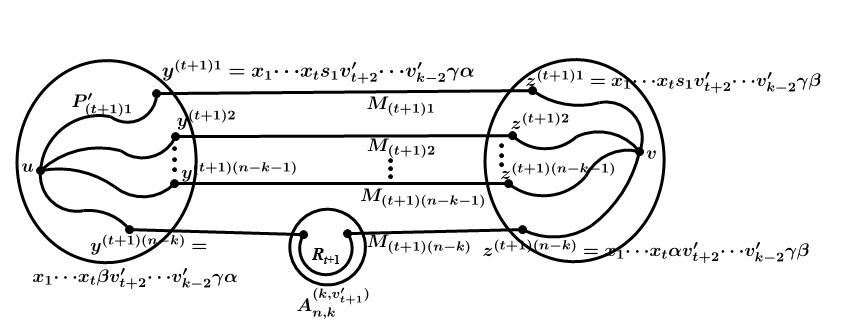}
  \includegraphics[width=0.88\textwidth]{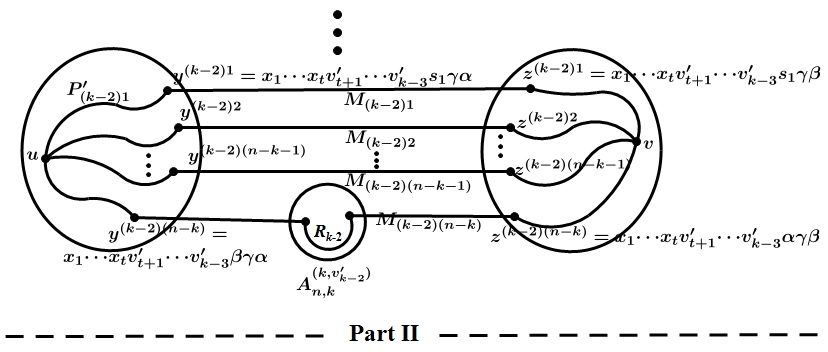}
 \end{figure}

(a) $(n-k)(k-1)+1\le l\le (n-k)(k-1)+n+t-2k+2$. (If $n+t-2k+2=0$, then (a) does not occur. )

\begin{figure}[!htbp]
 \centering
 \includegraphics[width=0.65\textwidth]{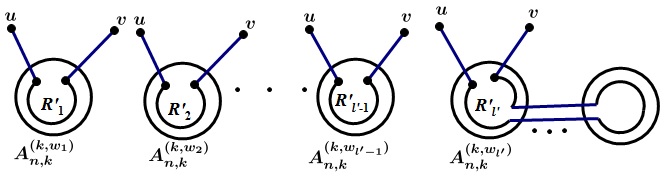}
 \caption{The paths $M_{k1}, \cdots, M_{kl'}$   of  case 2.1(a) in Lemma \ref{TH1} }
 \label{2_1_32}
 \end{figure}

 Let $l=(n-k)(k-1)+l'$, then $1\le l'\le n+t-2k+2$.
 By Lemma \ref{L2}, there exists a Hamiltonian path $R'_i$ of $A^{(k, w_i)}_{n, k}$ joining $u^{s(\alpha, w_i)}$ to $v^{s(\beta, w_i)}$ for $1\le i\le l'-1$. Let $I=\{ u'_{t+1}, \cdots, u'_{k-2}, w_{l'}, \cdots, w_{n+t-2k+2}\}$, by Lemma \ref{L2}, there exists a Hamiltonian path $R'_{l'}$ of $A^{(k, I)}_{n, k}$ joining $u^{s(\alpha, w_{l'})}$ to $v^{s(\beta, w_{l'})}$.
 Set $M_{ki}=\langle u, u^{s(\alpha, w_{i})}, R'_i, v^{s(\beta, w_i)}, v\rangle$ for $1\le i\le l'$.
 To construct the $l^*$-container, we only need to combine figure \ref{2_1_1} and $l'$ paths $M_{k1}, M_{k2}, \cdots, M_{kl'}$ in figure \ref{2_1_32}.

(b) $(n-k)(k-1)+n+t-2k+3\le l\le (n-k)k$. (If $k-t=2$, then (b) is not happened. )

Let $l=(n-k)(k-1)+n+t-2k+2+l'$, then $1\le l'\le k-t-2$.  In this case, we have $U\setminus(U\cap V)\not=\emptyset$.
Note that $\gamma\in \langle n\rangle\setminus V$. Without loss of generality, we can assure that $\gamma=u'_{k-2}$.

To construct the $l^*$-container,   we need to

\noindent{\bf Step 1: } By Lemma \ref{L2}, there exists a Hamiltonian path $R''_i$ of $A^{(k, w_i)}_{n, k}$ joining $u^{s(\alpha, w_i)}$ to $v^{s(\beta, w_i)}$ for $1\le i\le n+t-2k+2$.

\begin{figure}[!htbp]
 \centering
 \includegraphics[width=0.40\textwidth]{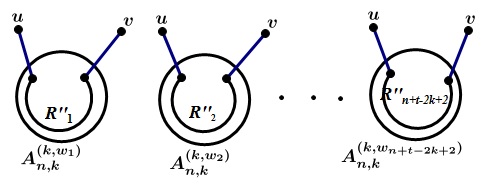}
 \caption{The paths $M_{k1}, \cdots, M_{k(n+t-2k+2)}$ of  case 2.1(b) in Lemma \ref{TH1}}
 \label{2_1_4}
 \end{figure}

  Combine figure \ref{2_1_1} and $(n+t-2k+2)$  paths $M'_{k1}, M'_{k2}, \cdots, M'_{k(n+t-2k+2)}$  in figure \ref{2_1_4}  where \begin{center}
$ M'_{ki}=  \langle u, u^{s(\alpha, w_i)}, R''_i, v^{s(\beta, w_i)}, v\rangle~~for~~1\le i\le n+t-2k+2 $.
\end{center}

 \noindent {\bf Step 2: } By Lemma \ref{L2}, there exists a Hamiltonian path $ H_{j}$
of  $ A^{(k, u'_{t+j})}_{n, k}$ joining $(y^{(t+j)(n-k)})^{s(\alpha, u'_{t+j})}$  to $ v^{s(\beta, u'_{t+j})} $  and a Hamiltonian path $ H'_{j}$ of $ A^{(k, v'_{t+j})}_{n, k}$ joining $u^{s(\alpha, v'_{t+j})}$ to $(z^{(t+j)(n-k)})^{s(\beta, v'_{t+j})}$ for $1\le j\le l'-1$.

\begin{figure}[!htbp]
 \centering
 \includegraphics[width=0.8\textwidth]{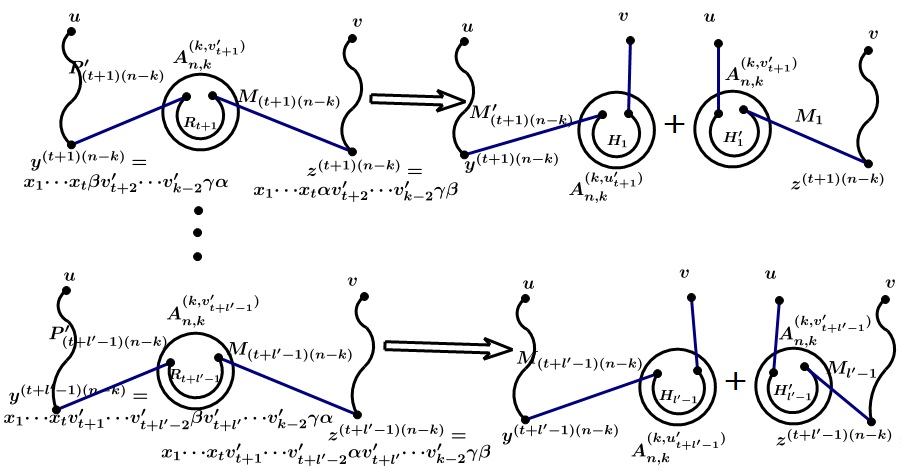}
 \caption{Illustration for  step 2 of  case 2.1(b) in Lemma \ref{TH1}}
 \label{2_1_5}
 \end{figure}

Replace paths $M_{(t+1)(n-k)}, M_{(t+2)(n-k)}, \cdots, M_{(t+l'-1)(n-k)}$ in part \Rmnum{2} of figure \ref{2_1_1}  by $M'_{(t+1)(n-k)}$, $M'_{(t+2)(n-k)}, \cdots$, $ M'_{(t+l'-1)(n-k)}$ and $M_1, M_2, \cdots$,$ M_{l'-1}$ as shown in figure \ref{2_1_5} where
    \begin{center}
      $\begin{array}{rl}
       M'_{(t+j)(n-k)}= & \langle u, P'_{(t+j)(n-k)}, y^{(t+j)(n-k)}, (y^{(t+j)(n-k)})^{s(\alpha, u'_{t+j})}, H_{j}, v^{s(\beta, u'_{t+j})}, v\rangle,\vspace{1.5ex}  \\
   M_j= & \langle u, u^{s(\alpha, v'_{t+j})}, H'_j, (z^{(t+j)(n-k)})^{s(\beta, v'_{t+j})}, z^{(t+j)(n-k)}, Q'_{(t+j)(n-k)}, v\rangle \vspace{1.5ex}
      \end{array}$
      for $1\le j\le l'-1$.
   \end{center}

 \noindent{\bf Step 3: } Let
$I=\{u'_{t+l'}, u'_{t+l'+1}, \cdots, u'_{k-2}\}$.
By Lemma \ref{L2}, there exists a Hamiltonian path $H$
of $A^{(k, u'_{k-2})}_{n, k}$ joining $u^{s(\alpha, v'_{k-2})}$ to $(z^{(k-2)(n-k)})^{s(\beta, v'_{k-2})}$ and a Hamiltonian path $ R$ of $ A^{(k, I)}_{n, k}$ joining $(y^{(k-1)(n-k)})^{s(\alpha, u'_{k-2})}$ to $v^{s(\beta, u'_{k-2})}$.

\begin{figure}[!htbp]
 \centering
 \includegraphics[width=0.8\textwidth]{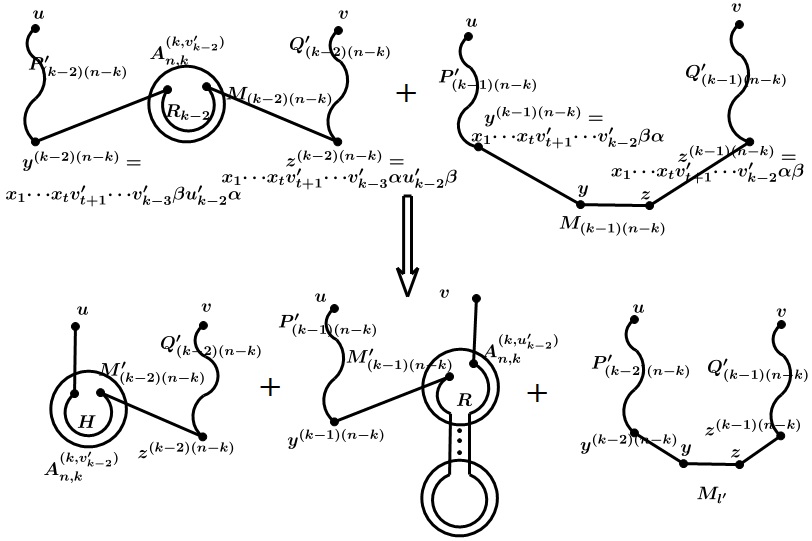}
 \caption{Illustration for  step 3 of  case 2.1(b) in Lemma \ref{TH1}}%{Illustration for case 2.1(b) of Theorem \ref{TH1}}
 \label{2_1_6}
 \end{figure}

 Replace $M_{(k-2)(n-k)}$ and $M_{(k-1)(n-k)}$ in  part \Rmnum{2} and part \Rmnum{3} of figure \ref{2_1_1}  by $M'_{(k-2)(n-k)}$, $M'_{(k-1)(n-k)}$  and $M_{l'}$
as shown  in figure \ref{2_1_6} where
\begin{center}
$\begin{array}{rl}
   M'_{(k-2)(n-k)}=& \langle u,  u^{s(\alpha, v'_{k-2})}, H, (z^{(k-2)(n-k)})^{s(\beta, v'_{k-2})}, z^{(k-2)(n-k)}, Q'_{(k-2)(n-k)}, v\rangle, \vspace{1.0ex}\\
  M'_{(k-1)(n-k)}= & \langle u, P'_{(k-1)(n-k)}, y^{(k-1)(n-k)}, (y^{(k-1)(n-k)})^{s(\alpha, u'_{k-2})}, R, v^{s(\beta, u'_{k-2})},v\rangle.\vspace{1.0ex}\\
  M_{l'}=&\langle u, P'_{(k-2)(n-k)}, y^{(k-2)(n-k)},  y, z, z^{(k-1)(n-k)}, Q'_{(k-1)(n-k)}, v\rangle.
 \end{array}
$
\end{center}
   \vspace{1.5ex}

\noindent{\bf Case 2.2 : } $\{i\mid (u)_i\notin V, (v)_i\in U\}\not=\emptyset$ or $(\{i\mid (v)_i\notin U, (u)_i\in V\}\not=\emptyset)$.

Without loss of generality, let $(u)_k\notin V, (v)_k\in U$,  and let $(u)_k=\alpha, (v)_k=\beta$.
Suppose that
\begin{center}
 $\begin{array}{rl}
   U\cap V= &\{  x_1, x_2, \cdots, x_t, \beta\}, \vspace{1.0ex}\\
     U\setminus (U\cap V)=& \{u'_{t+1}, u'_{t+2}, \cdots, u'_{k-2}, \alpha\}, \vspace{1.0ex} \\
      V\setminus (U\cap V)= &\{v'_{t+1}, v'_{t+2}, \cdots, v'_{k-1}\},\vspace{1.0ex} \\
      \langle n\rangle\setminus (U\cup V)=&\{w_1, w_2, \cdots, w_{n+t-2k+1}\} .
  \end{array}
 $
 \end{center}
 Without loss of generality, let $v=x_1x_2\cdots x_tv'_{t+1}\cdots v'_{k-1}\beta$.
 Since $n-k\ge 2$, there exists an element $\gamma\in\langle n\rangle\setminus V$. Set $y=x_1\cdots x_tv'_{t+1}\cdots v'_{k-2}\gamma \alpha$ and $z=x_1\cdots x_tv'_{t+1}\cdots v'_{k-2}\gamma \beta$. Thus, $u\not=y, z\not=v$.
Let
$S=\langle n\rangle\setminus\{\{(y)_i: 1\le i\le k\}\cup \{\beta\}\}=\langle n\rangle\setminus\{\{(z)_i: 1\le i\le k\}\cup\{\alpha\}\}=\{s_1, s_2, \cdots, s_{n-k-1}\}$.
By induction, there exists an $(n-k)(k-1)^*$-container $\{P_{ij}: 1\le i\le k-1, 1\le j\le n-k\}$
of $A^{(k, \alpha)}_{n, k}$ joining $u$ to $y$ and  an $(n-k)(k-1)^*$-container $\{Q_{ij}: 1\le i\le k-1, 1\le j\le n-k\}$
of $A^{(k, \beta)}_{n, k}$ joining $z$ to $v$.
We can represent
   $ P_{ij}$ as $ \langle u, P'_{ij}, y^{ij}, y\rangle, $
   $ Q_{ij}$ as $\langle z, z^{ij}, Q'_{ij}, v\rangle$
 for $1\le i\le k-1$ and $1\le j\le n-k$ where
\begin{center}
$ y^{ij}=\left\{
\begin{aligned}
  & y^{s((y)_i, s_j)}: 1\le i\le k-1, 1\le j\le n-k-1, \\
 &  y^{s((y)_i, \beta)}: 1\le i\le k-1, j= n-k,
\end{aligned}
\right.$
\end{center}
\begin{center}
$ z^{ij}=\left\{
\begin{aligned}
  & z^{s((y)_i, s_j)}: 1\le i\le k-1, 1\le j\le n-k-1, \\
 &  z^{s((y)_i, \alpha)}: 1\le i\le k-1,  j= n-k.
\end{aligned}
\right.$
\end{center}
Obviously,  $(y^{ij}, z^{ij})\in E(A_{n, k})$ when $1\le i\le k-1, 1\le j\le n-k-1$.
By Lemma \ref{L2}, there exists a Hamiltonian path  $R_{i}$ of $A^{(k, (y)_{i})}_{n, k}$ joining $(y^{i(n-k)})^{s(\alpha, (y)_i)}$ to $(z^{i(n-k)})^{s(\beta, (z)_i)}$ for $1\le i\le k-2$,
 Then, there exists $(n-k)(k-1)$ internally disjoint paths $\{M_{ij}: 1\le i\le k-1, 1\le j\le n-k\}$ of $A_{n, k}$ joining $u$ to $v$ such that
\begin{center}
$V(\displaystyle\bigcup_{i=1}^{k-1}\displaystyle\bigcup_{j=1}^{n-k} M_{ij})=V(A^{(k, \{\alpha, \beta, x_1, \cdots, x_t, v'_{t+1}, \cdots, v'_{k-2}\})}_{n, k})$.
\end{center}
where \begin{center}
$ M^{ij}=\left\{
\begin{aligned}
  &\langle u, P'_{ij}, y^{ij}, z^{ij}, Q'_{ij}, v\rangle~for~1\le i\le k-1, 1\le j\le n-k-1, \\
 &\langle u, P'_{ij}, y^{ij}, (y^{ij})^{s(\alpha, (y)_i)}, R_i, (z^{ij})^{s(\beta, (z)_i)}, z^{ij}, Q'_{ij}, v\rangle~for~1\le i\le k-2, j=n-k, \\
 &\langle u, P'_{ij}, y^{ij}, y, z, z^{ij}, Q'_{ij}, v\rangle~for~i=k-1, j=n-k.
\end{aligned}
\right.$
\end{center}
 See figure \ref{2_1_1} for illustration.

(a) $(n-k)(k-1)+1\le l\le (n-k)(k-1)+n+t-2k+2$.

Let $l=(n-k)(k-1)+l'$, then $1\le l'\le n+t-2k+2$.
 To construct the $l^*$-container, we need

\noindent{\bf Step 1: }By Lemma \ref{L2}, there exists a Hamiltonian path $R'_i$ of $A^{(k, w_i)}_{n, k}$ joining $u^{s(\alpha, w_i)}$ to $v^{s(\beta, w_i)}$ for $1\le i\le l'-1$.

\begin{figure}[!htbp]
 \centering
 \includegraphics[width=0.45\textwidth]{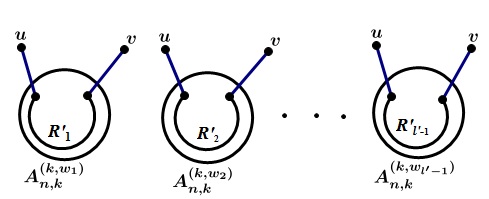}
 \caption{The paths $M_{k1}, \cdots, M_{k(l'-1)}$  of step 1 of  case 2.2(a) in Lemma \ref{TH1}}
 \label{2_2_3}
 \end{figure}

  Combine  figure \ref{2_1_1} and $(l'-1)$ disjoint paths $M_{k1}, M_{k2}, \cdots, M_{k(l'-1)}$  in figure \ref{2_2_3} where $M_{ki}= \langle u, u^{s(\alpha, w_i)}, R'_i, v^{s(\beta, w_i)}, v\rangle$ for $1\le i\le l'-1$.

\noindent{\bf Step 2:}
  Let $I=\langle n\rangle\setminus\{\alpha, \beta, x_1, \cdots, x_t, v'_{t+1}. \cdots, v'_{k-2}\}$. Then, there exists a Hamiltonian path $R'_{l'}$ of $A^{(k, I)}_{n, k}$ joining $u^{s(\alpha, v'_{k-1})}$ to $(z^{(k-1)(n-k)})^{s(\beta, v'_{k-1})}$.

\begin{figure}[!htbp]
 \centering
 \includegraphics[width=0.8\textwidth]{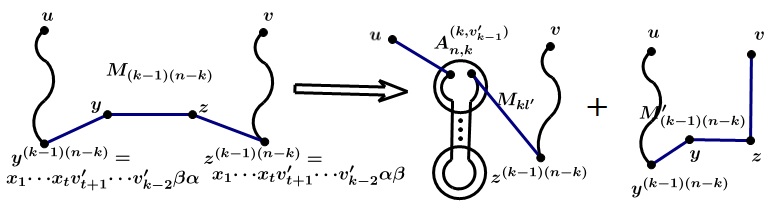}
 \caption{Illustration for  step 2 of  case 2.2(a) in Lemma \ref{TH1}}
 \label{2_2_4}
 \end{figure}

Replace $M_{(k-1)(n-k)}$ in  part \Rmnum{3} of figure \ref{2_1_1} by $M'_{(k-1)(n-k)}$ and $M_{kl'}$ as shown in \ref{2_2_4} where \vspace{1.0ex}\\
 $ \begin{array}{rl}
    M_{kl'}= & \langle u, u^{s(\alpha, v'_{k-1})}, R'_{l'}, (z^{(k-1)(n-k)})^{s(\beta, v'_{k-1})}, z^{(k-1)(n-k)}, Q'_{(k-1)(n-k)}, v\rangle, \vspace{1.0ex} \\
    M'_{(k-1)(n-k)}= & \langle u, P'_{(k-1)(n-k)}, y^{(k-1)(n-k)}, y, z, v\rangle\vspace{1.0ex}.
  \end{array}$

(b)$(n-k)(k-1)+n+t-2k+3\le l\le (n-k)(k-1)+n-k$. (If $k-t=2$, then (b) is not happened. )

Let $l=(n-k)(k-1)+n+t-2k+2+l'$, then $1\le l'\le k-t-2$. Thus, $U\setminus \{(U\cap V)\cup \{\alpha\}\}\not=\emptyset$,
 Note that $\gamma\in \langle n\rangle\setminus V$. Without loss of generality, we can assure that $\gamma=u'_{k-2}$.
 In order to construct the $l^*$-container, we need to

 \noindent{\bf Step 1: } By Lemma \ref{L2}, there exists a Hamiltonian path $R''_i$  of $A^{(k, w_i)}_{n, k}$ joining $u^{s(\alpha, w_i)}$
 to $v^{s(\beta, w_i)}$ for $1\le i\le n+t-2k+1$.

\begin{figure}[!htbp]
 \centering
 \includegraphics[width=0.48\textwidth]{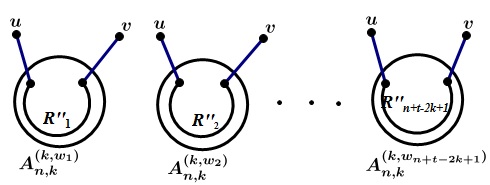}
 \caption{The paths $M'_{k1}, \cdots, M'_{k(n+t-2k+1)}$ of step 1 of  case 2.2(b) }
 \label{2_2_5}
 \end{figure}

  Combine figure \ref{2_1_1}  and $(n+t-2k+1)$ disjoint paths  $M'_{k1}, M'_{k2}, \cdots, M'_{k(n+t-2k+1)}$ in \ref{2_2_5} where $M'_{ki}=  \langle u, s^{s(\alpha, w_i)}, R''_{i}, v^{s(\beta, w_i)}, v\rangle$ for $1\le i\le n+t-2k+1\vspace{1.0ex}$.

 \noindent{\bf Step 2:} By Lemma \ref{L2},
 there exists a Hamiltonian path $H_{j}$ of  $A^{(k, u'_{t+j})}_{n, k}$
 joining $(y^{(t+j)(n-k)})^{s(\alpha, u'_{t+j})}$
 to  $v^{s(\beta, u'_{t+j})}$ and a Hamiltonian path $H'_{j}$ of $A^{(k, v'_{t+j})}_{n, k}$ joining $u^{s(\alpha, v'_{t+j})}$ to $(z^{(t+j)(n-k)})^{s(\beta, v'_{t+j})}$ for  $1\le j\le l'-1$.

\begin{figure}[!htbp]
 \centering
 \includegraphics[width=0.8\textwidth]{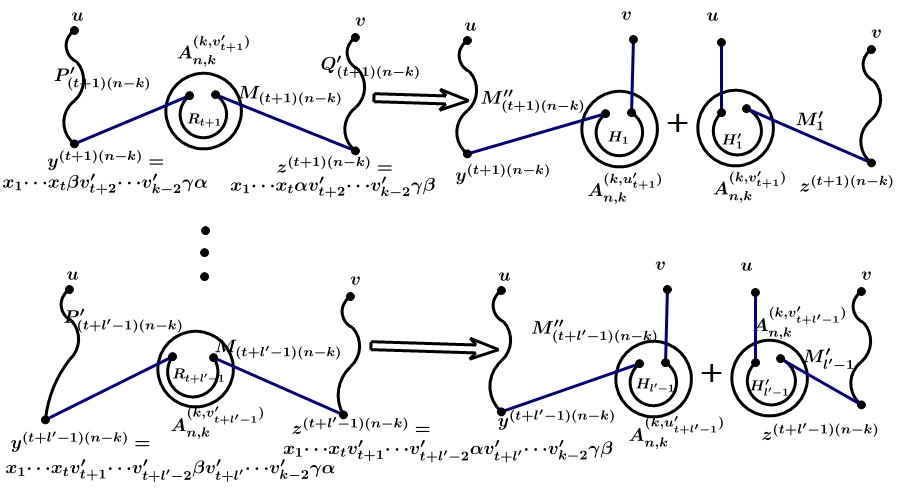}
 \caption{Illustration for  step 2 of  case 2.2(b) in Lemma \ref{TH1}}
 \label{2_2_6}
 \end{figure}

 Replace $M_{(t+1)(n-k)}, \cdots, M_{(t+l'-1)(n-k)}$ in  part \Rmnum{2} of figure \ref{2_1_1} by $M''_{(t+1)(n-k)}$, $\cdots,  M''_{(t+l'-1)(n-k)}$ and $M'_1, \cdots$, $M'_{l'-1}$ as shown in figure \ref{2_2_6} where
   \begin{center}
   $\begin{array}{rl}
      M'_{j}= & \langle u, u^{s(\alpha, v'_{t+j})}, H'_{j}, (z^{(t+j)(n-k)})^{s(\beta, v'_{t+j})},  z^{(t+j)(n-k)}, Q'_{(t+j)(n-k)}, v\rangle~for~1\le j\le l'-1, \vspace{1.0ex} \\
     M''_{(t+j)(n-k)}= & \langle u, P'_{(t+j)(n-k)}, y^{(t+j)(n-k)}, ( y^{(t+j)(n-k)})^{s(\alpha, u'_{t+j})} , H_{j}, v^{s(\beta, u'_{t+j})}, v\rangle ~for~1\le j\le l'-1.\\
   \end{array}$
   \end{center}

 \noindent{\bf Step 3: } By Lemma \ref{L2},
 there exists a Hamiltonian path  $ H_{k-1}$ of  $ A^{(k, u'_{k-2})}_{n, k}$
 joining  $(y^{(k-1)(n-k)})^{s(\alpha, u'_{k-2})}$
 to  $ v^{s(\beta, u'_{k-2})}$ and a Hamiltonian path $H'_{k-2}$ of $A^{(k, v'_{k-2})}_{n, k}$ joining $ u^{s(\alpha, v'_{k-2})}$ to $(z^{(k-2)(n-k)})^{s(\beta, v'_{k-2})}$.
  Let $I=\{v'_{k-1}, u'_{t+l'}, u'_{t+l'+1}, \cdots, u'_{k-2}\}$,  then there exists a Hamiltonian path  $H'_{k-1}$ of  $A^{(k, I)}_{n, k}$) joining  $ u^{s(\alpha, v'_{k-1})}$ to  $(z^{(k-1)(n-k)})^{s(\alpha, v'_{k-1})}$.

\begin{figure}[!htbp]
 \centering
 \includegraphics[width=0.8\textwidth]{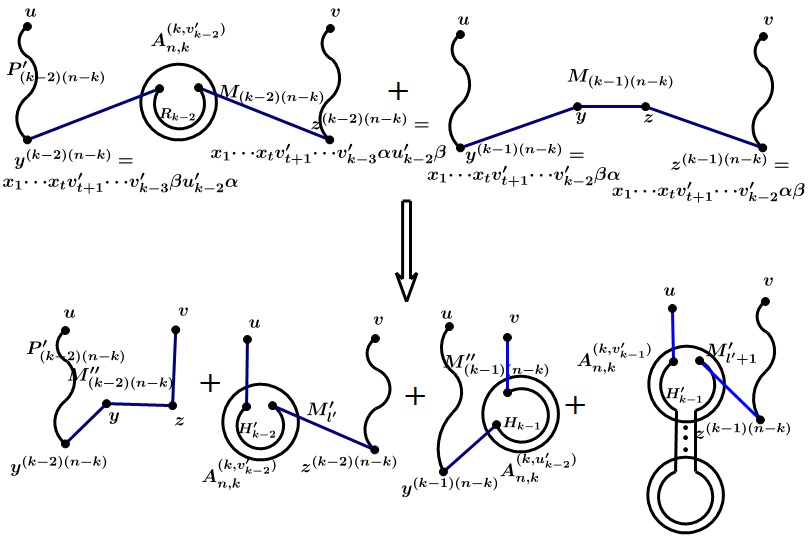}
 \caption{Illustration for  step 3 of  case 2.2(b) in Lemma \ref{TH1}}%{Illustration for case 2.1 of Theorem \ref{TH1}}
 \label{2_2_7}
 \end{figure}

 Replace $ M_{(k-2)(n-k)}, M_{(k-1)(n-k)}$ in  part \Rmnum{2} and part \Rmnum{3} of figure \ref{2_1_1}  by $ M''_{(k-2)(n-k)}$, $M''_{(k-1)(n-k)}$ and $M'_{l'}, M'_{l'+1}$ as shown in figure \ref{2_2_7} where
 \begin{flushleft}
   $\begin{array}{rl}
     M'_{l'}= &\langle u, u^{s(\alpha, v'_{k-2})}, H'_{k-2}, (z^{(k-2)(n-k)})^{s(\beta, v'_{k-2})}, z^{(k-2)(n-k)}, Q'_{(k-2)(n-k)}, v\rangle,\vspace{1.0ex}\\
     M'_{l'+1}=&\langle u, u^{s(\alpha, v'_{k-1})}, H'_{k-1}, (z^{(k-1)(n-k)})^{s(\beta, v'_{k-1})}, z^{(k-1)(n-k)}, Q'_{(k-1)(n-k)}, v\rangle,\vspace{1.0ex}\\
     M''_{(k-2)(n-k)}= & \langle u, P'_{(k-2)(n-k)}, y^{(k-2)(n-k)}, y, z, v\rangle, \vspace{1.0ex} \\
     M''_{(k-1)(n-k)}= &\langle u, P'_{(k-1)(n-k)}, y^{(k-1)(n-k)}, (y^{(k-1)(n-k)})^{s(\alpha, u'_{k-2})}, H_{k-1}, v^{s(\beta, u'_{k-2})}, v\rangle \vspace{1.0ex}.
   \end{array}$
 \end{flushleft}

\noindent{\bf Case 2.3 : } $U\cap V=\emptyset$.

Let $u=u_1u_2\cdots u_k, v=v_1v_2\cdots v_k, \langle n\rangle\setminus( U\cup V)=\{w_1, w_2, \cdots, w_{n-2k}\}$. Set $y=v_1v_2\cdots v_{k-1}u_k, z=u_1u_2\cdots u_{k-1}v_k$.
By induction, there exists an $(n-k)(k-1)^*$-container $\{P_{ij}: 1\le i\le k-1, 1\le j\le n-k\}$
of $A^{(k, u_k)}_{n, k}$ joining $u$ to $y$ and  an $(n-k)(k-1)^*$-container $\{Q_{ij}: 1\le i\le n-k, 1\le j\le n-k\}$
of $A^{(k, v_k)}_{n, k}$ joining $z$ to $v$.
We represent $P_{ij}$ as $\langle u, P'_{ij}, y^{ij}, y\rangle, Q_{ij}$ as $\langle z, z^{ij}, Q'_{ij}, v\rangle$ for $1\le i\le k-1, 1\le j\le n-k$ where
\begin{center}
  $ y^{ij}=\left\{
\begin{aligned}
  & y^{s(v_i, u_j)}: 1\le j\le k-1, \\
   & y^{s(v_i, v_k)}: j=k, \\
 &  y^{s(v_i, w_{j-k})}: k+1\le j\le n-k,
\end{aligned}
\right.
$
$ z^{ij}=\left\{
\begin{aligned}
 & z^{s(u_i, v_j)}: 1\le j\le k-1, \\
 & z^{s(u_i, u_k)}: j=k, \\
 &  z^{s(u_i, w_{j-k})}: k+1\le j\le n-k.
\end{aligned}
\right.
$
\end{center}

\noindent {\bf Subcase 2.3.1: }  $k=2$

Without loss of generality, let $u=12, v=34$.

When $n=5$, we set the $4^*$-container as :

 $~~~~~~~$ $\begin{array}{lll}
     \langle 12, 14, 34\rangle, ~~~~~~~~ &
   \langle12, 32, 34\rangle,  ~~~~~~~~&
  \langle12, 52, 34\rangle,   ~~~
   \end{array}\\
    $
    $~~~~~~~~~~~~~~~\langle12, 42, 43, 13, 53, 23, 21, 31, 51, 41$, $45, 35, 15, 25, 24, 54, 34\rangle.$

 Set the $5^*$-container as:

$~~~~~~~$ $\begin{array}{ll}
   \langle 12, 14, 34\rangle,& \langle 12, 15, 25, 45, 35, 34\rangle, ~~~~~~~  \langle 12, 32,  34\rangle, \\
    \langle 12, 52, 54,  34\rangle, &  \langle 12, 42, 43, 23, 13, 53, 51, 41, 31, 21, 24, 34\rangle.
 \end{array}$

Set the $6^*$-container as: \\
$~~~~~~~~~~~~~$ $\begin{array}{lll}
   \langle 12, 14, 34\rangle, &\langle 12, 13, 53, 43, 23, 24, 34\rangle, &  \langle 12, 15, 25, 45, 35, 34\rangle, \\
   \langle 12, 32, 34\rangle, & \langle 12, 42, 41, 51, 21, 31, 34\rangle, & \langle 12, 52, 54, 34\rangle.
 \end{array}
$\\
 Now, let $n\ge 6$.

(a) $l=(n-k)(k-1)+1=n-1$.

\begin{figure}[!htbp]
 \centering
 \includegraphics[width=0.45\textwidth]{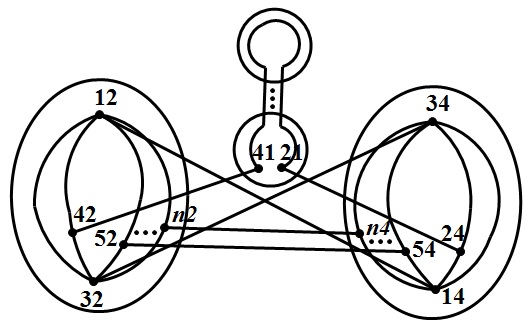}
 \caption{Illustration for subcase 2.3.1 (a) of Theorem \ref{TH1}}\label{2_3_1a}
 \end{figure}

 Let $I=n\setminus \{2, 4\}$, by Lemma \ref{L2}, there exists a Hamiltonian path $R_1$ of $A^{(2, I)}_{n, 2}$ joining $41$ to $21$, we set\\
 $\begin{array}{rlrl}
    P_1= &  \langle 12, 32, 34\rangle, &
   P_2= & \langle 12, 14, 34\rangle ,\\
    P_3=& \langle 12, 42, 41, R_1, 21, 24, 34\rangle,&
    P_{i}=&\langle 12, (i+1)2, (i+1)4, 34\rangle~~for~~ 4\le i\le n-1.
  \end{array}
 $\\
 See figure \ref{2_3_1a} for illustration.

(b) $l=(n-k)(k-1)+2=n$

\begin{figure}[!htbp]
 \centering
 \includegraphics[width=0.50\textwidth]{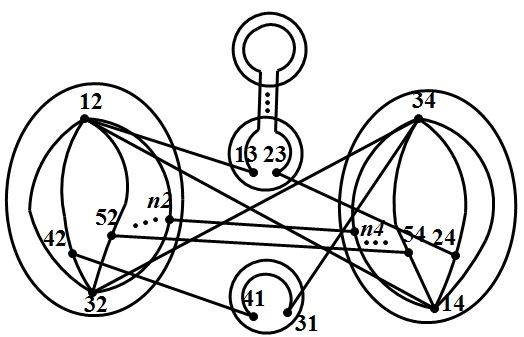}
 \caption{Illustration for subcase 2.3.1 (b) of Theorem \ref{TH1}}\label{2_3_1b}
 \end{figure}

By Lemma \ref{L2}, there exists a Hamiltonian path $R_1$ of $A^{(2, 1)}_{n, 2}$ joining $41$ to $31$. Let $I=\langle n\rangle\setminus\{1, 2, 4\}$, by Lemma \ref{L2}, there exists a Hamiltonian path $R_2$ of $A^{(2, I)}_{n, k}$ joining $13$ to $23$. We set\\
 $\begin{array}{rlrl}
    P_1= &  \langle 12, 32, 34\rangle,
   &P_2= & \langle 12, 14, 34\rangle ,\\
   P_3=& \langle 12, 42, 41, R_1, 31, 34\rangle,
   &P_4=&\langle 12, 13, R_2, 23, 24, 34\rangle,\\
    P_{i}=&\langle 12, i2, i4, 34\rangle~~for~~ 5\le i\le n.
  \end{array}
 $\\
  See figure \ref{2_3_1b} for illustration.

(c) $(n-k)(k-1)+3=n+1\le l\le (n-k)k= 2(n-2)$.

\begin{figure}[!htbp]
 \centering
 \includegraphics[width=0.55\textwidth]{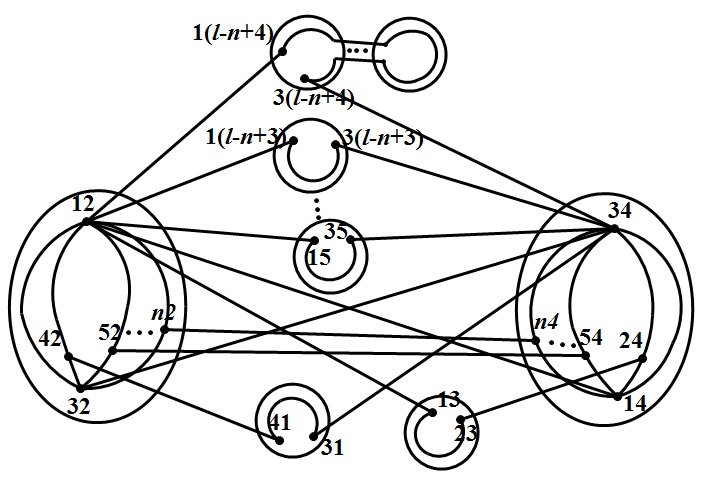}
 \caption{Illustration for subcase 2.3.1 (c) of Theorem \ref{TH1}}\label{2_3_1c}
 \end{figure}

By Lemma \ref{L2}, there exists a Hamiltonian path $R_1$ of $A^{(2, 1)}_{n, 2}$ joining  $41$ to $31$,  a Hamiltonian path $R_2$ of $A^{(2, 3)}_{n, k}$ joining $13$ to $23$,
and a Hamiltonian path $R_j$ of $A^{(2, j)}_{n, k}$ joining $1j$ to $3j$ for $5\le j\le l-n+3$.
Let $I=\langle n\rangle\setminus\{1, 2, \cdots, l-n+3\}$.
By Lemma \ref{L2}, there exists a Hamiltonian path $H_{l-n+4}$ of $A^{(2, l-n+4)}_{n, k}$ joining $1(l-n+4)$ to $3(l-n+4)$.
We set\begin{center}
 $\begin{array}{rl}
    P_1= &  \langle 12, 32, 34\rangle,\\
   P_2= & \langle 12, 14, 34\rangle ,\\
   P_3=& \langle 12, 42, 41, R_1, 31, 34\rangle,\\
   P_4=&\langle 12, 13, R_2, 23, 24, 34\rangle,\\
    P_{i}=&\langle 12, i2, i4, 34\rangle~~for~~ 5\le i\le n,\\
    P_{j}=&\langle 12, 1j, R_{j}, 3j, 34\rangle~~for~~ 5\le j\le l-n+3,\\
    P_{l-n+4}=&\langle 12, 1(l-n+4), H_{l-n+4}, 3(l-n+4), 34\rangle .
  \end{array}
 $\end{center}
  See figure \ref{2_3_1c} for illustration.

\noindent{\bf Subcase 2.3.2: } $k\ge 3$.

(a)\ $l=(n-k)(k-1)+1$.

Let
\begin{center}
 $\begin{array}{rl}
   A_1=&\{(y^{12})^{s(u_k, v_1)}, \cdots, (y^{1(k-1)})^{s(u_k, v_1)}, (y^{1(k+1)})^{s(u_k, v_1)},  \cdots, (y^{1(n-k)})^{s(u_k, v_1)}\}\vspace{1.0ex},\\
    B_1=&\{(z^{12})^{s(v_k, v_1)}, \cdots, (z^{1(k-1)})^{s(v_k, v_1)}, (z^{1(k+1)})^{s(v_k, v_1)},  \cdots, (z^{1(n-k)})^{s(v_k, v_1)}\} \vspace{1.0ex},\\
   A_i=&\{(y^{i1})^{s(u_k, v_i)}, \cdots, (y^{i(n-k)})^{s(u_k, v_i)}\}~~for~~2\le i\le k-1,\vspace{1.0ex} \\
   B_i=&\{(z^{i1})^{s(v_k, v_i)}, \cdots, (z^{i(i-1)})^{s(v_k, v_i)}, (z^{(i-1)(i-1)})^{s(v_k, v_i)}, (z^{i(i+1)})^{s(v_k, v_i)}, \cdots\vspace{1.0ex}, (z^{i(n-k)})^{s(v_k, v_i)}\}\\
   &~~~~~~~~~~~~~~~~~~~~~~~~~~~~~~~~~~~~~~~~~~~~~~~~~~~~~~~~~~for~~2\le i\le k-1. %\vspace{1.0ex}
  \end{array}
 $
\end{center}
We partite $A^{(k, v_i)}_{n, k}$ to $\cup_{j\in\langle n\rangle\setminus \{v_i\}}A^{(i, j)(k, v_i)}_{n, k}$.
By Lemma \ref{L6}, there exists $(n-k-2)$ disjoint paths $H_{12}, H_{13}, \cdots, H_{1(k-1)}, H_{1(k+1)}, \cdots, H_{1(n-k)}$ of $A^{(k, v_1)}_{n, k}$ from $A_1$ to $B_1$ such that  $V(\cup_{j\in \langle n-k\rangle\setminus\{1, k\}}H_{1j})=V(A^{(k, v_1)}_{n, k})$ and $H_{1j}=\langle (y^{1j})^{s(u_k, v_1)}, H_{1j}, (z^{1j})^{s(v_k, v_1)}\rangle$ for $j\in \langle n-k\rangle\setminus \{1, k\}$.
For $2\le i\le k-1$, there exists $(n-k)$ disjoint paths $H_{i1}, H_{i2}, \cdots, H_{i(n-k)}$ of $A^{(k, v_i)}_{n, k}$ from $A$ to $B$ such that $V(\cup_{j=1}^{n-k}H_{ij})=V(A^{(k, v_i)}_{n, k})$
and $H_{ii}=\langle (y^{ii})^{s(u_k, v_i)}, H_{ii}, (z^{(i-1)(i-1)})^{s(v_k, v_i)}\rangle$,
$H_{ij}=\langle (y^{ij})^{s(u_k, v_i)}, H_{ij}, (z^{ij})^{s(v_k, v_i)}\rangle$ for $j\in \langle n-k\rangle\setminus\{i\}$.
Let $I=\langle n\rangle \setminus\{v_1, \cdots, v_{k-1}, v_{k}$, $u_k\}$, by Lemma \ref{L2}, there exists a Hamiltonian path $H_{1k}$ of $A^{(k, I)}_{n,k}$ joining $(y^{1k})^{s(u_k, u_1)}$ to $(z^{1k})^{s(v_k, u_1)}$.  We set\vspace{1.0ex}

$\begin{array}{l}
   M_{11}=  \langle u, P'_{11}, y^{11}, y, v\rangle\vspace{1.0ex}, \\
   M_{1k}=  \langle u, P'_{1k}, y^{1k},  (y^{1k})^{s(u_k, u_1)}, H_{1k}, (z^{1k})^{s(v_k, u_1)}, z^{1k}, Q'_{1k}, v\rangle\vspace{1.0ex},\\
   M_{1t}=  \langle u, P'_{1t}, y^{1t}, (y^{1t})^{s(u_k, v_1)}, H_{1t}, (z^{1t})^{s(v_k, v_1)}, z^{1t}, Q'_{1t}, v\rangle~for~ 2\le t\le n-k~and~ t\not=k,\vspace{1.0ex} \\
   M_{ii}=  \langle u, P'_{ii}, y^{ii}, (y^{ii})^{s(u_k, v_i)}, H_{ii}, (z^{(i-1)(i-1)})^{s(v_k, v_i)}, z^{(i-1)(i-1)}, Q'_{(i-1)(i-1)}, v\rangle~~for~2\le i\le k-1, \vspace{1.0ex}\\
   M_{ij} = \langle u, P'_{ij}, y^{ij},  (y^{ij})^{s(u_k, v_i)}, H_{ij}, (z^{ij})^{s(v_k, v_i)}, z^{ij}, Q'_{ij}, v\rangle~for~2\le i\le k-1, \\
   ~~~~~~~~~~~~~~~~~~~~~~~~~~~~~~~~~~~~~~~~~~~~~~~~~~~~~~~~~~~~~~~~~~~~~~~~~1\le j\le n-k~and~j\not=i,\\
   M_{k1}=\langle u, z, z^{(k-1)(k-1)}, Q'_{(k-1)(k-1)}, v\rangle.\vspace{1.0ex}
 \end{array}$

\begin{figure}[!htbp]
 \centering
 \includegraphics[width=0.70\textwidth]{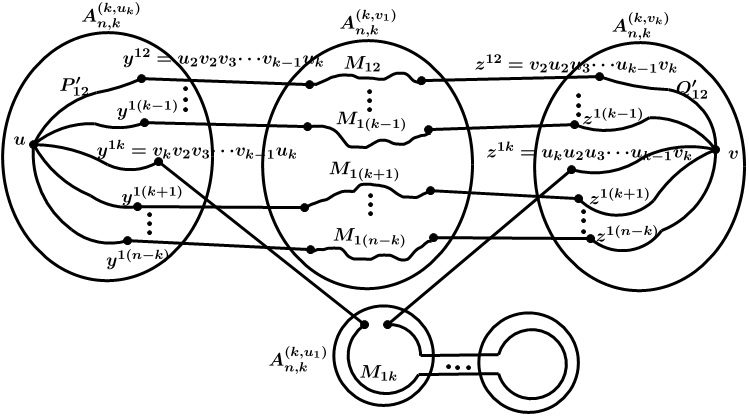}
  \includegraphics[width=0.70\textwidth]{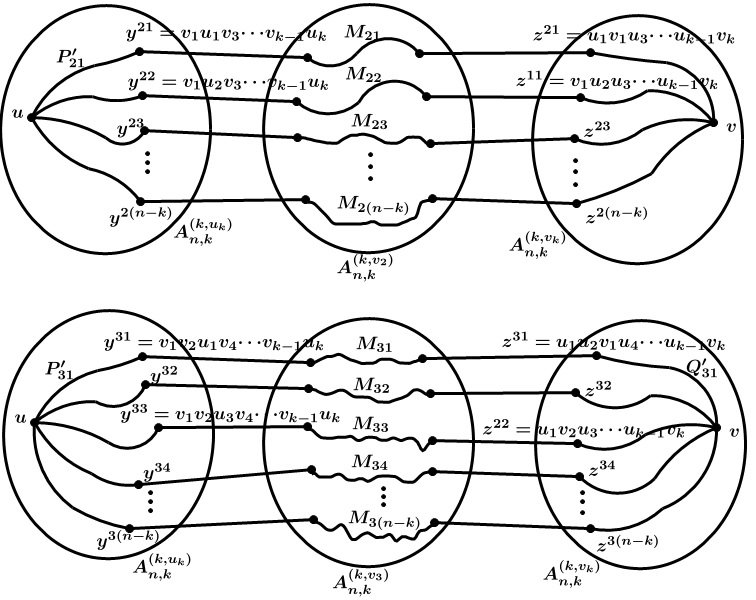}
   \includegraphics[width=0.80\textwidth]{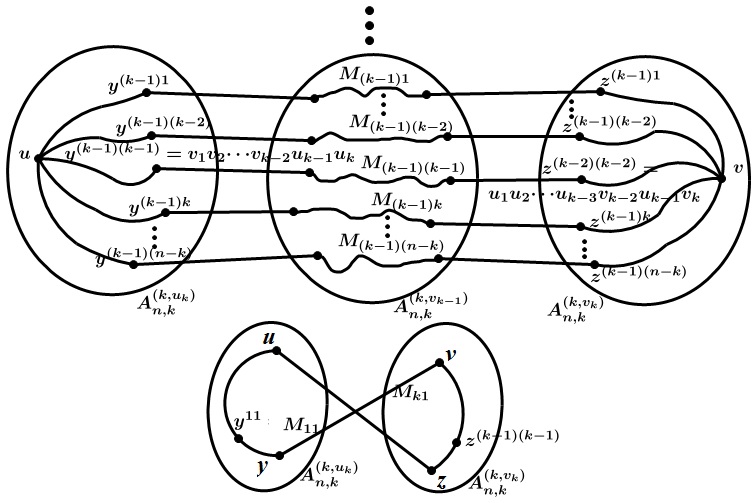}
 \caption{Illustration for subcase 2.3.2 (a) of Theorem \ref{TH1}}\label{2_3_1}
 \end{figure}

See figure \ref{2_3_1} for illustration.

(b) $(n-k)(k-1)+2\le l\le (n-k)(k-1)+n-2k+1$. (If $n-2k=0$, then (b) does not occur)

Let $l=(n-k)(k-1)+1+l'$ and
let
\begin{center}
 $\begin{array}{rl}
   A_1=&\{(y^{12})^{s(u_k, v_1)}, \cdots,  (y^{1(n-k)})^{s(u_k, v_1)}\}\vspace{1.0ex},\\
    B_1=&\{(z^{12})^{s(v_k, v_1)}, \cdots, (z^{1(n-k)})^{s(v_k, v_1)}\} \vspace{1.0ex}.
    \end{array}$\end{center}
%for $2\le j\le k-1$.
We partite $A^{(k, v_1)}_{n, k}$ to $\cup_{j\in\langle n\rangle\setminus \{v_1\}}A^{(1, j)(k, v_1)}_{n, k}$.
To construct the $l^*$-container, we need

\noindent{\bf Step 1: } By Lemma \ref{L6}, there exists $(n-k-1)$ disjoint paths $H'_{12}, H'_{13}, \cdots$, $H'_{1(n-k)}$ of $A^{(k, v_1)}_{n, k}$ from $A_1$ to $B_1$ such that $V(\cup_{j=2}^{n-k} H'_{1j})=V(A^{(k, v_1)}_{n, k})$ and $H'_{1j}=\langle (y^{1j})^{s(u_k, v_1)}, H'_{1j}, (z^{1j})^{s(v_k, v_1)}\rangle$ for $2\le j\le n-k$.

\begin{figure}[!htbp]
 \centering
 \includegraphics[width=0.70\textwidth]{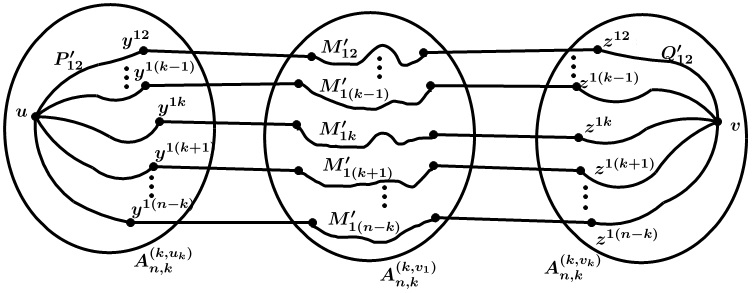}
 \caption{Illustration for  step 1 of  case 2.3(b) in Lemma \ref{TH1} }\label{2_3_4}
 \end{figure}

Replace paths $M_{12}, M_{13}, \cdots, M_{1(n-k)}$ in figure \ref{2_3_1} by $M'_{12}, M'_{13}, \cdots, M'_{1(n-k)}$ as shown in figure \ref{2_3_4} where
\begin{center}
$ M'_{1j}=  \langle u, P'_{1j}, y^{1j}, (y^{1j})^{s(u_k, v_1)}, H'_{1j}, (z^{1j})^{s(v_k, v_1)}, z^{1j}, Q'_{1j}, v\rangle~~for~~ j\in\langle n-k\rangle\setminus\{1\}.$
\end{center}

\begin{figure}[!htbp]
 \centering
 \includegraphics[width=0.70\textwidth]{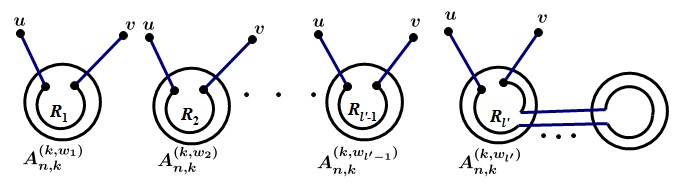}
 \caption{The paths $M_1, \cdots, M_{l'}$ of  step 2 of  case 2.3(b) in Lemma \ref{TH1} }\label{2_3_5}
 \end{figure}

\noindent{\bf Step 2: } Let $I=\langle n\rangle \setminus\{v_1, \cdots, v_{k-1}, v_{k}, u_k, w_1, \cdots, w_{l'-1}\}$.
By Lemma \ref{L2}, there exists a Hamiltonian path $R_{i}$ of $A^{(k, w_i)}_{n, k}$ joining $u^{s(u_k, w_{i})}$ to $v^{s(v_k, w_{i})}$ for $1\le i \le l'-1$
 and a Hamiltonian path $R_{l'}$ of $A^{(k, I)}_{n,k}$ joining $(u)^{s(u_k, w_{l'})}$ to $v^{s(v_k, w_{l'})}$.
Combine figure \ref{2_3_1} and  $l'$ paths $M_1, M_2, \cdots, M_{l'}$  in figure \ref{2_3_5} where
  \begin{center}
$
   M_{i}=\langle u, u^{s(u_k, w_{i})}, R_{i}, v^{s(v_k, w_{i})}, v\rangle$ for $ 1\le i\le l'.
$
\end{center}

(c)\ $(n-k)(k-1)+n-2k+2\le (n-k)(k-1)+n-k$.

Let $l=(n-k)(k-1)+n-2k+1+l'$, then $1\le l'\le k-1$.

To construct the $l^*$-container, we need:

\noindent{\bf Step 1: } By Lemma \ref{L2}, there exists a Hamiltonian path $R'_i$ of $A^{(k, w_i)}_{n, k}$ joining $u^{s(u_k, w_i)}$ to $v^{s(v_k, w_i)}$ for $1\le i\le n-2k$.

\begin{figure}[!htbp]
 \centering
 \includegraphics[width=0.50\textwidth]{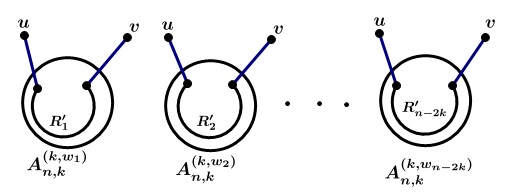}
 \caption{The paths $M'_1, \cdots, M'_{n-2k}$ of  step 1 of  case 2.3(c) in Lemma \ref{TH1}}\label{2_3_6}
 \end{figure}

 Combine figure \ref{2_3_1} and  $(n-2k)$ disjoint paths $M'_1, M'_{2}, \cdots, M'_{n-2k}$  as shown in figure \ref{2_3_6} where $M'_{i}=\langle u, u^{s(u_k, w_i)}, R'_i, v^{s(v_k, w_i)}, v\rangle$ for $1\le i\le n-2k$.
\bigskip

\bigskip

\noindent{\bf Step 2: } Notice $k\ge 3, n\ge 2k\ge 6$, then $\frac{(n-3)!}{(n-k-1)!}\ge 3$,
so there exists at least $3$ edges between $A^{(1, u_1)(k, v_1)}_{n, k}$ and $A^{(1, u_1)(k, u_{k-1})}_{n, k}$.
Let $(x^1, x^2)$ be an edge of $E(A^{(1, u_1)(k, v_1)}_{n, k}, A^{(1, u_1)(k, u_{k-1})}_{n, k})$ such that $x^1\in A^{(1, u_1)(k, v_1)}_{n, k}, x^2\in A^{(1, u_1)(k, u_{k-1})}_{n, k}$ and $x^1\not=u^{s(u_k, v_1)}, x^2\not=v^{s(v_k, u_{k-1})}$.
By Lemma \ref{L2}, there exists a Hamiltonian path $R$ of  $A^{(k, u_{k-1})}_{n, k}$ joining $x^2$ to $v^{s(v_k, u_{k-1})}$.
Let
\begin{center}
$\begin{array}{rl}
   A_1= & \{ u^{s(u_k, v_1)}, (y^{12})^{s(u_k, v_1)}, \cdots,  (y^{1(n-k)})^{s(u_k, v_1)}\}\vspace{1.0ex}, \\
   B_1= &\{ x^1, (z^{12})^{s(v_k, v_1)}, \cdots,  (z^{1(n-k)})^{s(v_k, v_1)}\}\vspace{1.0ex}.
 \end{array}
$\end{center}
By Lemma \ref{L6},
there exists $(n-k)$ disjoint paths $H''_{11}, H''_{12}, \cdots, H''_{1(n-k)}$ of $A^{(k, v_1)}_{n, k}$
from $A_1$ to $B_1$ such that $V(\cup_{i=1}^{n-k}H''_{1i})=V(A^{(k, v_1)}_{n, k})$ and $H''_{11}=\langle u^{s(u_k, v_1)}, H''_{11}, x^1\rangle,  H''_{1i}=\langle (y^{1i})^{s(u_k, v_1)}, H''_{1i}, (z^{1i})^{s(v_k, v_1)}\rangle$ for $2\le i\le n-k$.

\begin{figure}[!htbp]
 \centering
 \includegraphics[width=0.7\textwidth]{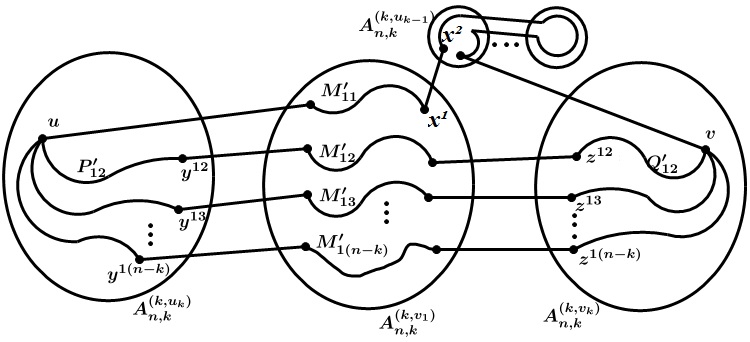}
 \caption{Illustration for step 2 of  case 2.3(c) in Lemma \ref{TH1}}\label{2_3_7}
 \end{figure}

Replace paths $M_{12}, M_{13}, M_{1(n-k)}$ in figure \ref{2_3_1} by paths $M''_{11}, M''_{12}, M''_{13}, \cdots, M''_{1(n-k)}$ as shown in figure \ref{2_3_7} where
\begin{center}
$\begin{array}{l}
  M''_{11}=  \langle u, u^{s(u_k, v_1)}, H''_{11}, x^1, x^2, R, v^{s(v_k, u_{k-1})}, v\rangle\vspace{1.0ex}, \\
   M''_{1i}=  \langle u, P'_{1i}, y^{1i}, (y^{1i})^{s(u_k, v_1)}, H''_{1i}, (z^{1i})^{s(v_k, v_1)}, z^{1i}, Q'_{1i}, v\rangle~for~ 2\le i\le n-k\vspace{1.0ex}.
   \end{array}$
\end{center}

\noindent{\bf Step 3: }Let
\begin{center}
$\begin{array}{rl}
  A_{i}= &\{(y^{i1})^{s(u_k, v_i)}, \cdots, (y^{i(i-1)})^{s(u_k, v_i)}, u^{s(u_k, v_i)}, (y^{i(i+1)})^{s(u_k, v_i)}, \cdots, (y^{i(n-k)})^{s(u_k, v_i)}\} \vspace{1.0ex} ~for~ 2\le i\le l'\\
   B^{i}= & \{(z^{i1})^{s(v_k, v_i)}, \cdots, (z^{i(i-1)})^{s(v_k, v_i)}, (z^{(i-1)(i-1)})^{s(v_k, v_i)}, (z^{i(i+1)})^{s(v_k, v_i)}, \cdots, (z^{i(n-k)})^{s(v_k, v_i)}\} \vspace{1.0ex}~\\
   &~~~~~~~~~~~~~~~~~~~~~~~~~~~~~~~~~~~~~~~~~~~~~~~~~~~~~~~~~~~~~~~~~~~~~~~~~~~~~~~~~~~ for~ 2\le i\le l'. \\
 \end{array}
$\\
\end{center}
For $2\le i\le l'$, we partite $A^{(k, v_i)}_{n, k}$ to $\cup_{r\in\langle n\rangle\setminus\{v_i\}}A^{(i, r)}_{n, k}$.  By Lemma \ref{L6}, there exists $(n-k)$ disjoint paths $H''_{i1}, H''_{i2}, \cdots, H''_{i(n-k)}$ of $A^{(k, v_i)}_{n, k}$ from $A^i$ to $B^i$ such that  $V(\cup_{j=1}^{n-k}H''_{ij})=V(A^{(k, v_i)}_{n, k})$ and
 $H''_{ii}=\langle u^{s(u_k, v_i)}, H''_{ii},  (z^{(i-1)(i-1)})^{s(v_k, v_i)}\rangle $,
 $H''_{ij}=$$\langle (y^{ij})^{s(u_k, v_i)}, H''_{ij}, (z^{ij})^{s(v_k, v_i)}\rangle$ for $j\in \langle n-k\rangle\setminus \{i\}$.
 By Lemma \ref{L2}, there exists a Hamiltonian path $H_{ki}$ of $A^{(k, u_{i-1})}_{n, k}$ joining $(y^{ii})^{s(u_k, u_{i-1})}$ to $v^{s(v_k, u_{i-1})}$ for $2\le i\le l'$.

\begin{figure}[!htbp]
 \centering
 \includegraphics[width=0.8\textwidth]{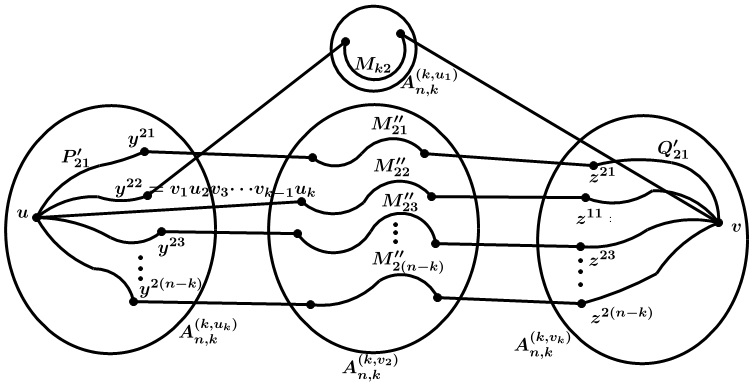}
 \includegraphics[width=0.8\textwidth]{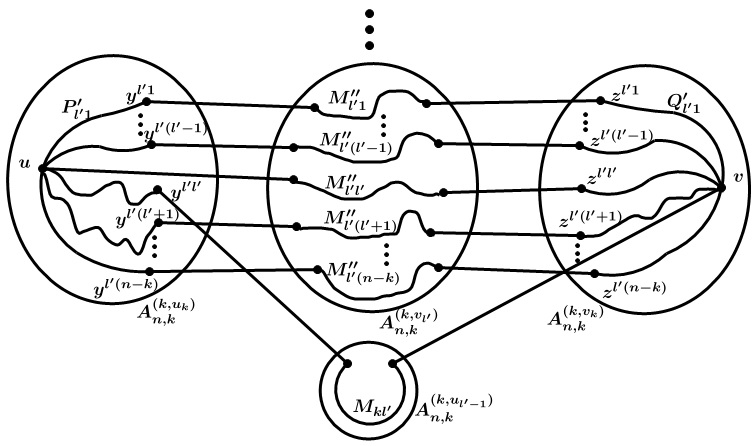}
 \caption{Illustration for  step 3 of  case 2.3(c) in Lemma \ref{TH1}}\label{2_3_8}
 \end{figure}

 Replace paths $M_{i1}, M_{i2}, \cdots, M_{i(n-k)}$ in figure \ref{2_3_1} by $M_{ki}, M''_{i1}, M''_{i2}, \cdots, M''_{i(n-k)}$ for $2\le i\le l'$ as shown in figure \ref{2_3_8}  where
\begin{center}
$\begin{array}{rl}
   M_{ki}= & \langle u, P'_{ii}, y^{ii}, (y^{ii})^{s(u_k, u_{i-1})}, H_{ki}, v^{s(v_k, u_{i-1})}, v\rangle, \vspace{1.0ex}\\
  M''_{ii}=& \langle u, u^{s(u_k, v_i)}, H''_{ii}, (z^{(i-1)(i-1)})^{s(v_k, v_i)}, z^{(i-1)(i-1)}, Q'_{(i-1)(i-1)}, v\rangle, \vspace{1.0ex}\\
  M''_{ij}=&\langle u, P'_{ij}, y^{ij}, (y^{ij})^{s(u_k, v_i)}, H''_{ij}, (z^{ij})^{s(v_k, v_i)}, z^{ij}, Q'_{ij}, v\rangle.\vspace{1.0ex}
 \end{array}
$ \\ for $2\le i\le l', 1\le j\le n-k~and~j\not=i. $
\end{center}

$\qed$

\begin{thm}{\rm }\label{TH2}
   $A_{n, k}$ is super spanning connected for $n\ge 4, n-k\ge2$. %$l^*$-connected for $3\le l\le(n-k)(k-1)$.
\end{thm}
{\bf Proof}: We prove the theorem by induction.

\noindent{\bf Basis step:} It is known that $A_{n, 1}$ is isomorphic to the complete graph $K_n$ and by Lemma \ref{L5}, $A_{4,2}$ is super spanning connected. Then, the result holds for $A_{n, 1}$ and $A_{4, 2}$.

\noindent{\bf Induction step:} Suppose $A_{n-1, k-1}$ is super spanning connected. We need to prove that $A_{n, k}$ is super spanning connected for $n\ge 5, n-k\ge2$. By Lemma \ref{L1}, $A_{n, k}$ is $1^*$-connected and $2^*$-connected.
By Lemma \ref{TH1}, $A_{n, k}$ is $l^*$-connected for $(n-k)(k-1)+1\le l\le (n-k)k$.
Now, we need to construct a $l^*$-container of $A_{n, k}$ joining any two distinct vertices $u$ and $v$ for $3\le l\le (n-k)(k-1)$. We use $U$ to denote the  set $\{(u)_i\mid 1\le i \le k\}$ and use $V$ to denote the set $\{(v)_i\mid 1\le i\le k\}$.

\noindent{\bf Case 1:} $\{i\mid (u)_i=(v)_i: 1\le i\le k\}\not=\emptyset$.

\begin{figure}[!htbp]
 \centering
 \includegraphics[width=0.40\textwidth]{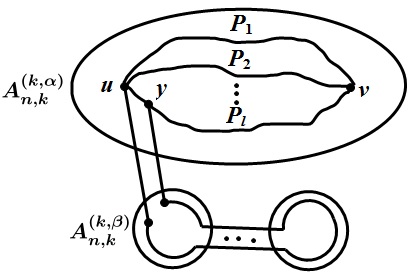}
 \caption{Illustration for case 1 of Theorem \ref{TH2}}\label{TH2case1}
 \end{figure}

Without loss of generality, let $(u)_k=(v)_k=\alpha$. By induction, there is an $l^*$-container $\{P_1, P_2, \cdots, P_l\}$ of
 $A^{(k, \alpha)}_{n, k}$ joining $u$ to $v$. Hence, we can represent $P_l$ as $\langle u, y, P'_l, v\rangle$.
 Note that $|\{(u)_i: 1\le i\le k \}\cup \{(y)_i: 1\le i\le k \}|=k+1$ and $n-k\ge 2$. Suppose  $\beta\in\langle n\rangle\setminus\{\{(u)_i: 1\le i\le k \}\cup \{(y)_i: 1\le i\le k \}\}$.
 By Lemma \ref{L2}, there exists a Hamiltonian path $H$ of $A^{(k, \langle n\rangle\setminus\{\alpha\})}_{n, k}$ joining $u^{s(\alpha, \beta)}$ to $y^{s(\alpha, \beta)}$.
 We set $P''_l=\langle u, u^{s(\alpha, \beta)}, H, y^{s(\alpha, \beta)}, y, P'_l, v\rangle$.
 Obviously, $\{P_1, P_2, \cdots, P_{l-1}, P''_l\}$ is  a $l^*$-container of $A_{n, k}$ joining $u$ to $v$. See figure \ref{TH2case1} for illustration.

\noindent{\bf Case 2:} $\{i\mid (u)_i=(v)_i: 1\le i\le k\}=\emptyset$.

\noindent{\bf Case 2.1:} $U\not=V$.

Without loss of generality, we can assume that $(u)_k=\alpha\notin V$.  We partite $A_{n, k}$ to $\cup_{i\in \langle n\rangle}A^{(k, i)}_{n, k}$. Suppose $u=u_1u_2\cdots u_{k-1}\alpha, v=v_1v_2\cdots v_{k-1}\beta$. Set $y=v_1v_2\cdots v_{k-1}\alpha$, then $y\not=u$ and $(y, v)\in E(A_{n, k})$. By induction, there exists an $l^*$-container $\{P_1, P_2, \cdots, P_l\}$ of $A^{(k, \alpha)}_{n, k}$ joining $u$ to $y$. We represent $P_i$ as $\langle u, P'_i, y^i, y\rangle$ for $1\le i\le l$. Without loss of generality , we can assume that $V(P_l)\le V(P_i)$ for $1\le i\le l$.
Suppose $|\{y^i\mid \beta\in\{(y^i)_j: 1\le j\le k-1\}, 1\le i\le l-1\}|=m$, then $0\le m\le$ min $\{l-1, k-1\}$.

\noindent{\bf Subcase 2.1.1} $m=0$.

\begin{figure}[!htbp]
 \centering
 \includegraphics[width=0.65\textwidth]{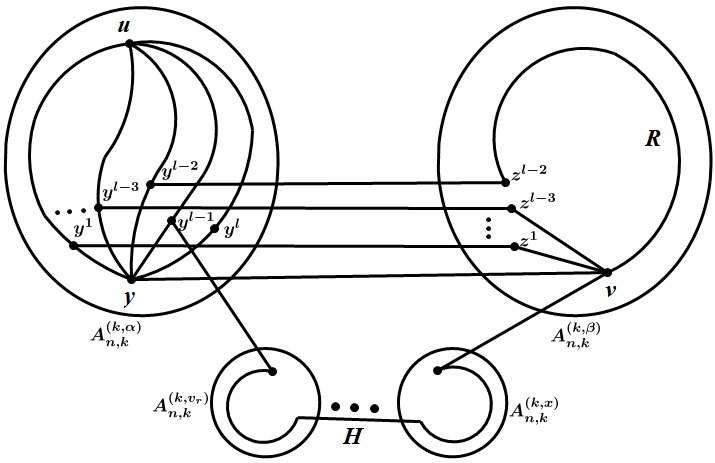}
 \caption{Illustration for subcase 2.1.1 of Theorem \ref{TH2}}\label{3_2_1}
 \end{figure}

For all $1\le i\le l-1$, we have $\beta\notin\{(y^i)_j: 1\le j\le k-1\}$. Let $z^i=(y^i)^{s(\alpha, \beta)}$ for $1\le i\le l-2$. Then, $z^i\in A^{(k, \beta)}_{n,k}$ and $(z^i, v)\in E(A_{n, k})$ for $1\le i\le l-2$. Note that $l-3\le (n-k)(k-1)-3$, by Lemma \ref{L1}, there exists a Hamiltonian path $R$ of $A^{(k, \beta)}_{n, k}\setminus\{z^1, z^2, \cdots, z^{l-3}\}$ joining $z^{l-2}$ to $v$. Since $y^{l-1}$ and $y$ differ in exactly one position,  we can assume that $y^{l-1}=v_1\cdots v_{r-1} x v_{r+1}\cdots v_{k-1}\alpha$ where $x\in \langle n\rangle\setminus\{v_1, v_2, \cdots, v_{k-1}, \alpha, \beta\}$.
Let $I=\langle n\rangle\setminus\{\alpha, \beta\}$. By Lemma \ref{L2}, there exists a Hamiltonian path $H$ of
$A^{(k, I)}_{n, k}$ joining $(y^{l-1})^{s(\alpha, v_r)}$ to $v^{s(\beta, x)}$. We set\\
  $\begin{array}{rl}
     M_i= & \langle u, P'_i, y^i, z^i, v\rangle ~for~1\le i\le l-3, \\
     M_{l-2}=& \langle u, P'_{l-2}, y^{l-2}, z^{l-2}, R, v\rangle,  \\
     M_{l-1}= & \langle u, P'_{l-1}, y^{l-1}, (y^{l-1})^{s(\alpha, v_r)}, H, v^{s(\beta, x)}, v\rangle, \\
     M_{l}= & \langle u, P'_l, y^{l}, y, v\rangle.
   \end{array}
  $\\
Obviously, $\{M_1, M_2, \cdots, M_l\}$ forms an $l^*$-container of $A_{n, k}$. See figure \ref{3_2_1} for illustration. 

\noindent{\bf Subcase 2.1.2} $0<m<l-1$.

\begin{figure}[!htbp]
 \centering
 \includegraphics[width=0.70\textwidth]{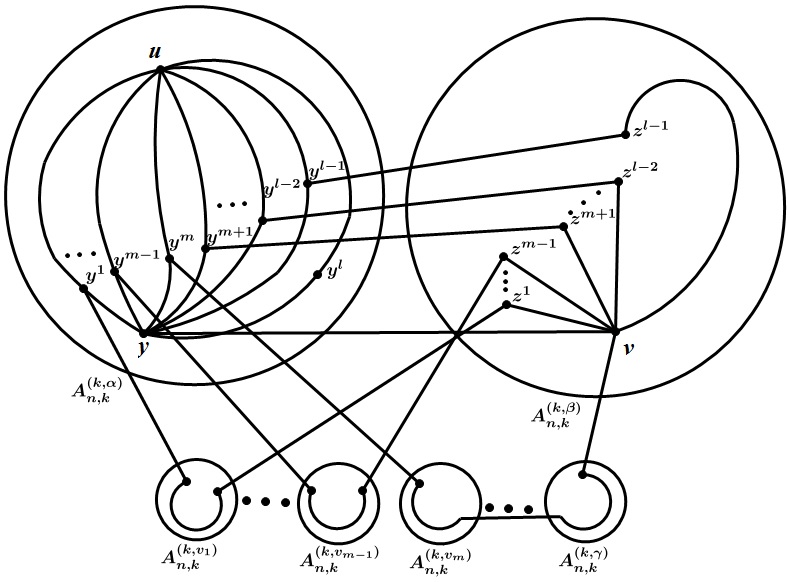}
 \caption{Illustration for subcase 2.1.2 of Theorem \ref{TH2}}\label{3_2_2}
 \end{figure}

Without loss of generality, we can assume that $y^i=v_1\cdots v_{i-1}\beta v_{i+1}\cdots v_{k-1}\alpha$ for $1\le i\le m$. Set $z^{i}=v_1\cdots v_{i-1}\alpha v_{i+1}\cdots v_{k-1}\beta$ for $1\le i\le m$. By Lemma \ref{L2}, there exists a Hamiltonian path $R_i$ of $A^{(k, v_i)}_{n, k}$ joining $(y^i)^{s(\alpha, v_i)}$ to $(z^i)^{s(\beta, v_i)}$ for $1\le i\le m-1$. Since $n-k\ge 2$, suppose $\gamma\in\langle n\rangle\setminus\{v_1, \cdots, v_{k-1}, \alpha, \beta\}$. Let $I=\langle n\rangle\setminus\{v_1, \cdots, v_{m-1}, \alpha, \beta\}$. By Lemma \ref{L2}, there exists a Hamiltonian path $R_m$ of $A^{(k, I)}_{n, k}$ joining $(y^{m})^{s(\alpha, v_m)}$ to $v^{s(\beta, \gamma)}$.
Notice that $\beta\notin\{(y^i)_j: 1\le j\le k-1\}$ for $m+1\le i\le l-1$.
We set $z^i=(y^i)^{s(\alpha, \beta)}$ for $m+1\le i\le l-1$. Then, $(y^i, z^i)\in E(A_{n, k})$ for $m+1\le i\le l-1$.
Since $l-3\le (n-k)(k-1)-3$, by Lemma \ref{L1}, there exists a Hamiltonian path $R$ of $A^{(k, \beta)}_{n, k}\setminus\{v_1, \cdots, v_{m-1}, v_{m+1}, \cdots, v_{l-2}\}$ joining $z^{l-1}$ to $v$. We set\\
$
\begin{array}{rl}
  M_{i}= & \langle u, P'_i, y^i, (y^i)^{s(\alpha, v_i)}, R_i, (z^{i})^{s(\beta, v_i)}, z^i, v\rangle~for~1\le i\le m-1,\vspace{1.0ex} \\
  M_{m}= & \langle u, P'_m, y^m, (y^{m})^{s(\alpha, v_m)}, R_m, v^{s(\beta, \gamma)}, v\rangle,\vspace{1.0ex} \\
  M_{m+j}= & \langle u, P'_{m+j}, y^{m+j}, z^{m+j}, z\rangle~for~1\le j\le l-m-2,\vspace{1.0ex} \\
  M_{l-1}= & \langle u, P'_{l-1}, y^{l-1}, z^{l-1}, R, v\rangle,\vspace{1.0ex} \\
  M_l=& \langle u, P'_{l}, y^l, y, v\rangle.
\end{array}
$\\
Obviously, $\{M_1, M_2, \cdots, M_l\}$ forms an $l^*$-container of $A_{n, k}$. See figure \ref{3_2_2} for illustration. 

\noindent{\bf Subcase 2.1.3} $m=l-1$.

 Without loss of generality, we can assume $y^i=v_1\cdots v_{i-1}\beta v_{i+1}\cdots v_{k-1}\alpha$ for $1\le i\le l-1$. Set $z^{i}=v_1\cdots v_{i-1}\alpha v_{i+1}\cdots v_{k-1}\beta$ for $1\le i\le l-1$. By Lemma \ref{L2}, there exists a Hamiltonian path $R_i$ of $A^{(k, v_i)}_{n, k}$ joining $(y^i)^{s(\alpha, v_i)}$ to $(z^i)^{s(\beta, v_i)}$ for $1\le i\le l-2$. Since $n-k\ge 2$, suppose $\gamma\in\langle n\rangle\setminus\{v_1, \cdots, v_{k-1}, \alpha, \beta\}$. Let $I=\langle n\rangle\setminus\{v_1, \cdots, v_{l-2}, \alpha, \beta\}$. By Lemma \ref{L2}, there exists a Hamiltonian path $R_m$ of $A^{(k, I)}_{n, k}$ joining $(y^{l-1})^{s(\alpha, v_{l-1})}$ to $(z^{l-1})^{s(\beta, \gamma)}$. Since $1\le l-2=m-1\le k-2\le(n-k)(k-1)-3$, by Lemma \ref{L2}, there exists a Hamiltonian path $R$ of $A^{(k, \beta)}_{n, k}$ joining $z^{l-1}$ to $v$. We set

\begin{figure}[!htbp]
 \centering
 \includegraphics[width=0.7\textwidth]{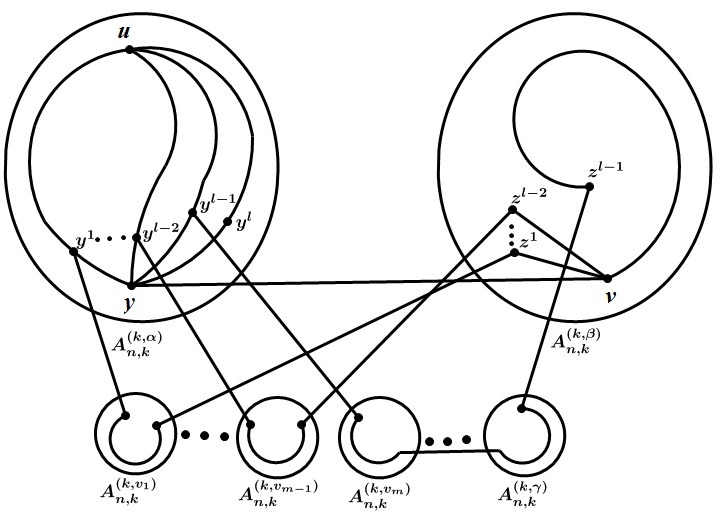}
 \caption{Illustration for subcase 2.1.3 of Theorem \ref{TH2}}\label{3_2_3}
 \end{figure}

$
\begin{array}{rl}
  M_i= & \langle u, P'_i, y^i, (y^i)^{s(\alpha, v_i)}, R_i, (z^i)^{s(\beta, v_i)}, z^i, v\rangle ~for~1\le i\le l-2,\vspace{1.0ex} \\
  M_{l-1}= & \langle u, P'_{l-1}, y^{l-1}, (y^{l-1})^{s(\alpha, v_{l-1})}, R_{l-1}, (z^{l-1})^{s(\beta, v_{l-1})}, z^{l-1}, R, v\rangle, \vspace{1.0ex} \\
  M_{l}=& \langle u, P'_l, y^l, y, v\rangle.
\end{array}
$\\
Obviously, $\{M_1, M_2, \cdots, M_l\}$ forms an $l^*$-container of $A_{n, k}$. See figure \ref{3_2_3} for illustration. 

 \noindent{\bf Case 2.2} $U=V$.

 \noindent{\bf Subcase 2.2.1} $k=2$.

\begin{figure}[!htbp]
 \centering
 \includegraphics[width=0.50\textwidth]{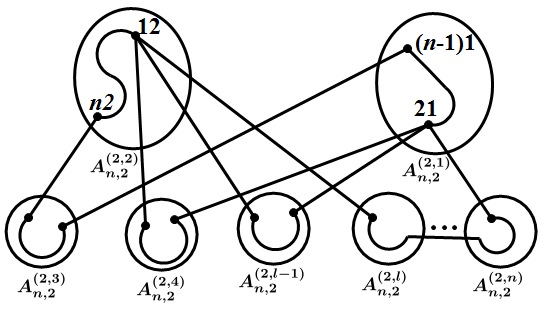}
 \caption{Illustration for subcase 2.2.1 of Theorem \ref{TH2}}\label{3_2_4}
 \end{figure}

 Without loss of generality, we can assume that $u=12, v=21, n\ge 5$ . By Lemma \ref{L2}, there exists a Hamiltonian path $R_1$ of $A^{(2, 2)}_{n, 2}$ joining $12$ to $n2$, a Hamiltonian path $R_2$ of $A^{(2, 1)}_{n, 2}$ joining $(n-1)1$ to $21$ and a Hamiltonian path $R_3$ of $A^{(2, 3)}_{n, k}$ joining $n3$ to $(n-1)3$. Additionally, there exists a Hamiltonian path $R_i$ of $A^{(2, i)}_{n, 2}$ joining $1i$ to $2i$ for $4\le i\le l-1$ and a Hamiltonian path $R$ of $A^{(k, \{l, l+1, \cdots, n\})}_{n, 2}$ joining $1l$ to $2n$.
  We set \\
  $\begin{array}{rl}
     M_1= & \langle 12, R_1, n2, n3, R_3, (n-1)3, (n-1)1, R_2, 21\rangle, \\
     M_i= &\langle 12, 1(i+2), R_{i+2}, 2(i+2), 21\rangle~ for ~2\le i\le l-1,\\
     M_{l}=&\langle 12, 1l, R, 2n, 21\rangle.
   \end{array}
  $\\
 Obviously, $\{M_1, M_2, \cdots, M_l\}$ forms an $l^*$-container of $A_{n, 2}$. See figure \ref{3_2_4} for illustration. 
 \bigskip

 \noindent{\bf Subcase2.2.2} $k\ge 3$.

 Without loss of generality, we can assume that $v=12\cdots k$, $u=u_1u_2\cdots u_{k-1}r$ and $r=1$. Let $y=k23\cdots(k-1)1$.
   By induction, there exists an $l^*$-container $\{P_1, P_2, \cdots, P_{l}\}$ joining $u$ to $y$. We represent $P_i$ as $\langle u, P'_i, y^i, y\rangle$ for $1\le i\le l$.

 If $|\{y^i\mid (y^i)_1\not=k\}|=1$, we may assume that $(y^l)_1\not=k$. If $|\{y^i\mid (y^i)_1\not=k\}|=2$, we may assume that $(y^{l-1})_1\not=k, (y^l)_1\not=k$.

 %Let $Y^j=\{y^i\mid (y^i)_j\not=(y)_j: 1\le i\le l-2\}=\{y^{j1}, y^{j2}, \cdots, y^{jn_j}\}$ for $1\le j\le k-1$.
 Let $Y=Y^1\cup Y^2\cup \cdots \cup Y^{k-1}$ be the neighbors of $y$ in $P_i$ for $1\le i\le l-2$ where $Y^j=\{y^{j1}, y^{j2}, \cdots, y^{j{n_j}}\}$ are  obtained by switch the $j$th coordinate of $y$.
 Then, $l-2=n_1+n_2+\cdots+n_{k-1}$. Now, we use
 $~\{P_{11}, P_{12}, \cdots, P_{1n_1}, P_{21}, P_{22}, \cdots, P_{2n_2}, \cdots, P_{(k-1)1}, P_{(k-1)2}, \cdots, P_{(k-1)n_{k-1}}, P_{l-1}, P_{l}\}$ to denote the $l^*$-container where $P_{ij}=\langle u, P'_{ij},  y^{ij}, y\rangle$ for $i\in\langle k-1\rangle, j\in\langle n_i\rangle$ and $P_{l-1}=\langle u, P'_{l-1}, y^{l-1}, y\rangle$, $ P_{l}=\langle u, P'_l, y^l, y\rangle$. We use $z^{ij}$ to denote the vertex $z^{s(i, (y^{ij})_i)}$ for $i=2, \cdots, k-1$( for example: if $y^{ij}=k2\cdots (i-1)x(i+1)\cdots(k-1)1$, then $z^{ij}=12\cdots(i-1)x(i+1)\cdots k$ ).
Thus, $(z^{ij}, v)\in E(A_{n, k})$.

 For $i\in \langle k\rangle\setminus\{1, r\}$, if $n_i\not=0$, we partite $A^{(k, i)}_{n, k}$ to $\cup_{j\in \langle n\rangle\setminus\{i\}}A^{(i, j)(k, i)}_{n, k}$. Let
 \begin{center}
  $ \begin{array}{cc}
      A^i= & \{(y^{i1})^{s(1, i)}, (y^{i2})^{s(1, i)}, \cdots, (y^{in_i})^{s(1, i)}\}, \vspace{1.0ex} \\
      B^i= & \{(z^{i1})^{s(k, i)}, (z^{i2})^{s(k, i)}, \cdots, (z^{in_i})^{s(k, i)}\}
    \end{array}$
 \end{center}
 By Lemma \ref{L6}, there exists $n_i$ disjoint paths $H_{i1}, H_{i2}, \cdots, H_{in_i}$ from $A^i$ to $B^i$ such that
 \begin{center}
 $V(\displaystyle\bigcup_{j=1}^{n_i}H_{ij})=V(A^{(k, i)}_{n, k})$ and $H_{ij}=\langle (y^{ij})^{s(1, i)}, H_{ij}, (z^{ij})^{s(k, i)}\rangle$ for $1\le j\le n_j$.
 \end{center}

 (a) $\{y^i\mid (y^i)_1\not=(y)_1\}\le 2$.

 Note that $l-2=n_1+n_2+\cdots n_{k-1}\ge 1$ and $n_1=0$. Without loss of generality, we can assume that $n_{k-1}\not=0$. Since $l-3\le (n-k)(k-1)-3$, by Lemma \ref{L1}, there exists a Hamiltonian path $R$ of $A^{(k, k)}_{n, k}\setminus \{z^{21}, \cdots, z^{2n_2}, \cdots, z^{(k-1)1}, \cdots, z^{(k-1)(n_{k-1}-1)}\}$
joining $z^{(k-1)n_{k-1}}$ to $v$. We set \vspace{1.0ex}

\begin{center}
$M_{ij}=\left\{
\begin{aligned}
  & \langle u, P'_{ij},  y^{ij}, (y^{ij})^{s(1, i)}, H_{ij}, (z^{ij})^{s(k, i)}, z^{ij}, v \rangle: i=2, 3, \cdots, k-2, j=1, 2, \cdots, n_i , \\
 & ~~~~~~~~~~~~~~~~~~~~~~~~~~~~~~~~~~~~~~~~~~~~~~~~~~~or~~ i=k-1, 1\le j\le n_{k-1}-1,\\
 &  \langle u, P'_{ij},  y^{ij}, (y^{ij})^{s(1, i)}, H_{ij}, (z^{ij})^{s(k, i)}, z^{ij}, R, v \rangle: i=k-1, j=n_{k-1}.
\end{aligned}
\right.$
\end{center}

Since $n-k\ge2$, let $a\in \langle n\rangle\setminus\{\{(y)_i: 1\le i\le k\}\cup\{(y^{l-1})_i: 1\le i\le k\}\}, b\in \langle n\rangle\setminus\{1, 2, \cdots, k, a \}$.
Let $I=\langle n\rangle\setminus(\{i\mid n_i\not=0: 2\le i\le k-1\}\cup\{k, a\})$. By Lemma \ref{L2}, there exists a Hamiltonian path $R_{l-1}$ of $A^{(k, a)}_{n, k}$ joining $(y^{l-1})^{s(1, a)}$ to $v^{s(k, a)}$ and a Hamiltonian path $R_l$ of $A^{(k, I)}_{n, k}$ joining $y^{s(1, b)}$ to $v^{s(k, b)}$.
We set

\begin{figure}[!htbp]
 \centering
 \includegraphics[width=0.75\textwidth]{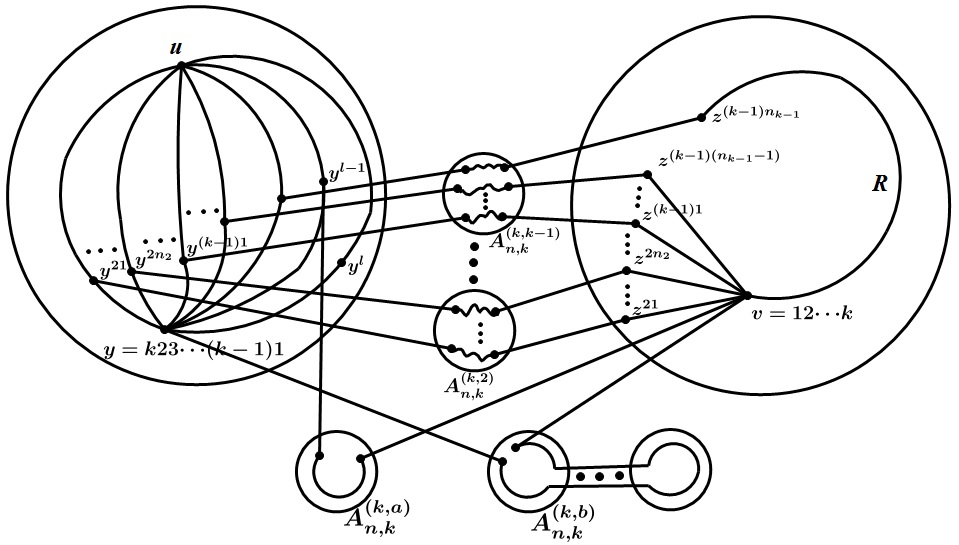}
 \caption{Illustration for subcase 2.2.2(a) of Theorem \ref{TH2}}\label{3_2_5}
 \end{figure}

$M_{l-1}=\langle u, P'_{l-1}, y^{l-1}, (y^{l-1})^{s(1, a)}, R_{l-1}, v^{s(k, a)}, v\rangle\vspace{1.0ex}$,

$M_l=\langle u, P'_l ,y^l, y, y^{s(1, b)}, R_{l}, v^{s(k, b)}, v\rangle$. \\
Obviously, $\{M_{21}, \cdots, M_{2n_2}, \cdots, M_{(k-1)1}, \cdots, M_{(k-1)}n_{k-1}, M_{l-1}, M_{l}\}$ forms an $l^*$-container of $A_{n, k}$. See figure \ref{3_2_5} for illustration. 

(b) $\{y^i\mid (y^i)_1\not=(y)_1\}\ge 3$.

Note that $n_2+\cdots n_{k-1}\le l-3\le (n-k)(k-1)-3$ .  By Lemma \ref{L1}, there exists a Hamiltonian path $R$ of $A^{(k, k)}_{n, k}\setminus \{z^{21}, \cdots, z^{2n_2}, \cdots, z^{(k-1)1}, \cdots, z^{(k-1)n_{k-1}}\}$
joining $z^{11}$ to $v$. Let $X=\langle n\rangle\setminus \langle k\rangle=\{x_1, x_2, \cdots, x_{n-k}\}$, Without loss of generality,  we may assume that $y^{1i}=y^{s(k, x_i)}$ for $1\le i\le n_1$, $y^{l-1}=y^{s(k, x_{n_1+1})}, y^l=y^{s(k, x_{n_1+2})}$ and $z^{11}=v^{s(1, x_{n_1+1})}$.
By Lemma \ref{L2}, there exists a Hamiltonian path $R_i$ of $A^{(k, x_{i+1})}_{n, k}$ joining $(y^{1i})^{s(1, x_{i+1})}$ to $v^{s(k, x_{i+1})}$ for $1\le i\le n_1$ and a Hamiltonian path $R_{l-1}$ of $A^{(k, x_{n_1+2})}_{n, k}$ joining $(y^{l-1})^{s(1, x_{n_1+2})}$ to $v^{s(k, x_{n_1+2})}$.
Let $I=\langle n\rangle\setminus\{\langle k\rangle\cup \{x_2, x_3,  \cdots, x_{n_1+2}\}\}$. By Lemma \ref{L2}, there exists a Hamiltonian $\vspace{1.0ex}$ path $R_l$ of  $A^{(k, I)}_{n, k}$ joining $(y)^{s(1, x_1)}$ to $(z^{11})^{s(k, x_1)}$. We set\vspace{1.0ex}

\begin{figure}[!htbp]
 \centering
 \includegraphics[width=0.8\textwidth]{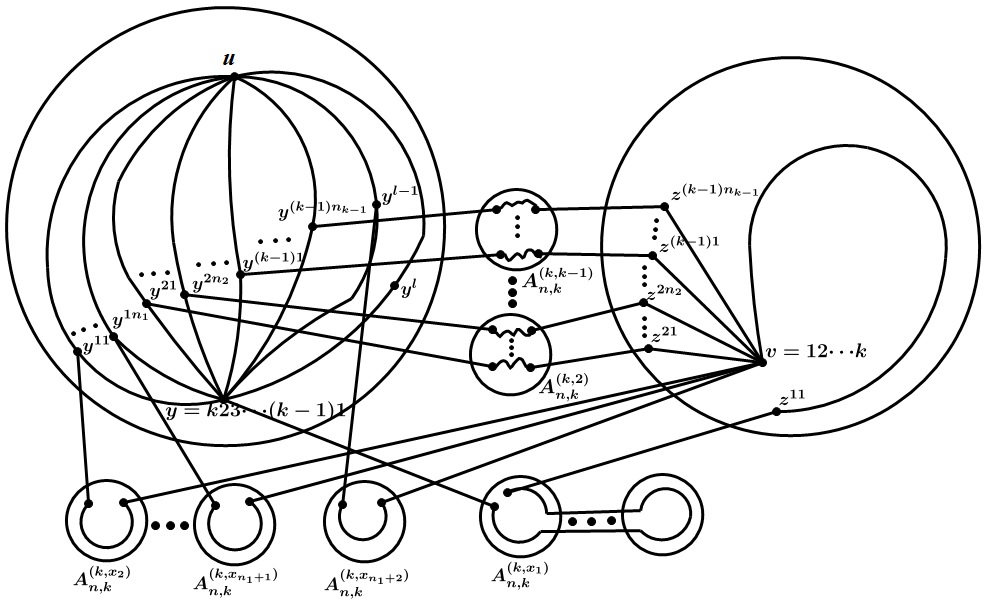}
 \caption{Illustration for subcase 2.2.2(b) of Theorem \ref{TH2}}\label{3_2_6}
 \end{figure}

$\begin{array}{rl}
  M_{1t}=&\langle u, P'_{1t}, y^{1t}, (y^{1t})^{s(1, x_{t+1})}, R_t, v^{s(k, x_{t+1})}, v\rangle$ for $1\le t\le n_1\vspace{1.0ex}, \\
  M_{ij}=&
   \langle u, P'_{ij},  y^{ij}, (y^{ij})^{s(1, i)}, H_{ij}, (z^{ij})^{s(k, i)}, z^{ij}, v \rangle~for~ 2\le i\le k-1, 1\le j\le n_i.\\
M_{l-1}=&\langle u, P'_{l-1}, y^{l-1}, (y^{l-1})^{s(1, x_{n_1+2})}, R_{l-1}, v^{s(k, x_{n_1+2})}, v\rangle\vspace{1.0ex},\\
M_l=&\langle u, y^l, y, y^{s(1, x_1)}, R_{l}, (z^{11})^{s(k, x_1)}, z^{11}, R, v\rangle.
 \end{array}
$\\
Obviously, $\{M_{11}, \cdots, M_{1n_1}, M_{21}, \cdots, M_{2n_2}, \cdots, M_{(k-1)1}, \cdots, M_{(k-1)}n_{k-1}, M_{l-1}, M_{l}\}$ forms an $l^*$-container of $A_{n, k}$. See figure \ref{3_2_6} for illustration. 
$\qed$

\section*{Acknowledgement}
This research is supported by National Natural Science Foundation of China (11571044, 61373021),  the Fundamental Research Funds for the Central University of China and Guangxi Provincial Education Department(YB2014007).

\end{CJK*}

\end{document}